\documentclass[preprint]{elsarticle}
\usepackage[utf8]{inputenc}

\usepackage{float}
\usepackage[margin=1in]{geometry}
\usepackage{comment}
\usepackage{stackengine}
\usepackage{lmodern}

\usepackage[thinlines]{easytable}
\usepackage{algorithm}
\usepackage{algpseudocode}
\usepackage{svg}
\usepackage{subcaption}
\usepackage{graphicx}
\usepackage{epsfig}
\usepackage{epstopdf}

\usepackage{amssymb}

\usepackage{bm}
\usepackage{amsmath}
\usepackage{amsfonts}
\usepackage{physics}
\usepackage{xfrac}
\usepackage[bb = dsserif]{mathalpha}
\usepackage[mathscr]{euscript}
\renewcommand{\comment}[2][.2\linewidth]{%
  \leavevmode\hfill\makebox[#1][l]{~#2}}

\usepackage[]{bm}

\usepackage[breaklinks, colorlinks=false,linkcolor=blue, citecolor=blue, urlcolor=blue]{hyperref}
\usepackage[capitalise]{cleveref} 
\usepackage{doi} 
\biboptions{sort}
\usepackage[page,toc]{appendix} 
\crefname{appsec}{Appendix}{Appendix}

\begin{document}
\begin{frontmatter}

\title{Spacetime Wavelet Method for the Solution of Nonlinear Partial Differential Equations}
\author{Cody D. Cochran}
\author{Karel Matou\v{s}\corref{cor1}}
\cortext[cor1]{Corresponding author}
\journal{Journal}

\address{Department of Aerospace and Mechanical Engineering, University of Notre Dame, Notre Dame, IN 46556, USA}
\begin{abstract}
\noindent
We propose a high-order spacetime wavelet method for the solution of
nonlinear partial differential equations with a user-prescribed
accuracy. The technique utilizes wavelet theory with \textit{a priori}
error estimates to discretize the problem in both the spatial and
temporal dimensions simultaneously. We also propose a novel
wavelet-based recursive algorithm to reduce the system sensitivity
stemming from steep initial and/or boundary conditions.  The resulting
nonlinear equations are solved using the Newton-Raphson method. We
parallelize the construction of the tangent operator along with the
solution of the system of algebraic  equations. We perform rigorous
verification studies using the nonlinear Burgers' equation. The
application of the method is demonstrated solving Sod shock tube
problem using the Navier-Stokes equations. The numerical results of
the method reveal high-order convergence rates for the function as
well as its spatial and temporal derivatives. We solve problems with
steep gradients in both the spatial and temporal directions with
\textit{a priori} error estimates. 
\end{abstract}


\end{frontmatter}
\section{Introduction} \label{sec:intro}
This work outlines the mathematical formulation, computational
implementation, and numerical results of the high-order spacetime
wavelet method as it applies to the solution of nonlinear partial
differential equations (PDEs), modeling shock behavior. Traditional
computational methods are capable of solving various problems, however
computational resources such as computing power and memory continue to
be limiting factors \cite{ferreau2008online,
  torabzadehkashi2019computational}. Simulations of nonlinear systems
involving multiple spatial and temporal scales and steep gradients,
such as those caused by shocks, tend to encounter one or both of these
limitations. This can result in long run times, insufficient accuracy,
or an inability to obtain predictive results
\cite{woodward1984numerical, sundarapandian2018strong,
  figuli2020numerical}.  

Many numerical methods have been developed to spatially discretize and
solve systems of nonlinear PDEs, including finite element methods
(FEM) \cite{gui1985h, dong2003p}, multigrid methods
\cite{thekale2010optimizing, dendy1982black, yushu2020image}, and
adaptive mesh refinement (AMR) techniques \cite{berger1984adaptive,
  klein1999star}, for example. These techniques possess several
advantages, such as their ability to adapt the grid to complex
geometry (\textit{e.g.}, FEM) and an optimal number of operations that
scales linearly with the number of nodes (\textit{e.g.},
multigrid). However, these approaches are computationally expensive
(\textit{e.g.}, AMR) in the context of problems with multiple spatial
and temporal scales where the location of steep solution features is
not known \textit{a priori}. 

This work develops an alternative discretization strategy, based on wavelets, to resolve all relevant scales and provide high-order convergence. In their early development as mathematical tools, wavelets were used primarily for image compression and signal processing \cite{daubechies1992ten}. These applications opened the door to exploring other areas for which wavelets can be employed. In particular, wavelet-based methods have been utilized to solve PDE systems \cite{bacry1992wavelet, goedecker2009wavelets, beylkin1997adaptive, vasilyev2000second, bertoluzza1996wavelet, harnish2018adaptive, harnish2021multiresolution, harnish2023adaptive}.

Wavelet approximations use a series of weighted basis functions similar to spectral methods \cite{bernardi1997spectral}. The wavelet theory is based on Multiresolution Analysis (MRA) which provides rigorous mathematical structure, including \textit{a priori} error estimates for an approximation of a function and its derivatives \cite{goedecker2009wavelets}. Wavelet-based PDE solvers have been developed and have merit in multiscale modeling \cite{schneider2010wavelet}, data compression \cite{beylkin1997adaptive}, nonlinear PDE solutions \cite{paolucci2014wamr1, paolucci2014wamr2, harnish2023adaptive}, bounded energy conservation \cite{ueno2003wavelet, qian1993wavelets}, and error control \cite{harnish2018adaptive, harnish2023adaptive}. As mentioned above, wavelets provide a consistent approximation of functions with \textit{a priori} error estimates. However, computing derivatives of functions is a more complex process, and most wavelet-based solvers use methods such as the finite difference method to provide an approximation of the derivative terms \cite{paolucci2014wamr1, paolucci2014wamr2, alam2006simultaneous}. Although finite difference approximations used properly in conjunction with wavelet analysis have been shown to provide adequate solutions to PDEs, discretizing derivative terms with wavelets provides a global high-order approximation as predicted by wavelet theory \cite{harnish2018adaptive, harnish2023adaptive}. 

Typical time-dependent PDE solution techniques involve semi-discretization into a system of ordinary differential equations (ODEs) \cite{moin2010fundamentals}, leaving a decision to be made about how to integrate in time. When solving PDEs, the time-integration process can become a major source of error and a computational bottleneck if not implemented appropriately. Most solution techniques integrate in time using explicit or implicit methods \cite{dokainish1989survey, subbaraj1989survey}, which can become inefficient when dealing with a stiff system of equations \cite{loffeld2013comparative}.  The timestep interval required to keep an explicit scheme stable can be so small that it becomes impractical to use \cite{ashi2008numerical}. Implicit schemes tend to accept larger time steps and can be unconditionally stable for certain, mostly linear, problems, but also require more computational effort \cite{loffeld2013comparative}. Several implicit, parallel time-integration techniques have been developed that involve a decomposition of the time domain, where each subdomain is solved asynchronously \cite{gander2007analysis, gander201550, subber2016asynchronous, prakash2014computationally, benevs2010asynchronous}. These time integrators have been shown to obtain accurate solutions for a variety of problems with truncation error estimates expressed as functions of the time step, \textit{i.e.}, $\mathcal{O}\left(\Delta t ^n\right)$. Despite having a predictable truncation error, other sources can cause the global error to deviate from local estimates \cite{verwer1996explicit, loffeld2013comparative}.

Discretization methods exist that do not require time marching techniques and are capable of computing the full solution in both space and time in a single solve. These methods are often referred to in the literature as \textit{spacetime} methods \cite{abedi2004spacetime, hughes1996space}. With spacetime methods, the spatial and temporal dimensions are discretized similarly and simultaneously, removing the need for time-stepping. Published work demonstrates the validity of spacetime methods as useful PDE solvers. Many of these methods use variations of FEM \cite{hulbert1990space, abedi2010adaptive, takizawa2012space, hughes1996space, petersen2009space}, but some work has also been done utilizing wavelets \cite{alam2006simultaneous, gunzburger2011space, schwab2009space}. However, like many other wavelet PDE solvers, time derivatives are usually approximated with more conventional techniques (\textit{e.g.}, finite difference), forfeiting \textit{a priori} global error estimates provided by wavelet theory. 

This work introduces a novel spacetime wavelet method that uses wavelet basis functions and wavelet derivatives to fully formulate and discretize systems of nonlinear PDEs in both space and time. This removes the need for semi-discretization and time-marching. The residual strong form system of nonlinear equations is solved using the Newton-Raphson method. In contrast to time-marching schemes, we do not have a stability condition that dictates the timestep size. However, we introduce a solvability condition that assesses the sensitivity of the system of linear algebraic equations and controls the construction of the initial guess and the tangent matrix conditioning. In this regard, we introduce a wavelet-based recursive algorithm that progressively builds the solution from a less steep initial and/or boundary condition and uses wavelet synthesis to construct a viable initial guess. We demonstrate the high-order convergence rates not only for the function but also for all spatial and temporal derivatives in the PDE solution. Furthermore, we control user-prescribed accuracy of the solution utilizing \textit{a priori} error estimates. 

The remainder of this paper is organized as follows: Section \ref{sec:wavelet} breaks down the mathematical formulation of wavelet theory and its application to discretize differential equations. Next, Section \ref{sec:nonlinearComp} describes the numerical and computational techniques used to implement and solve spacetime wavelet problems. Lastly, Section \ref{sec:numerical} presents and analyzes the numerical results of the spacetime wavelet technique as applied to a variety of nonlinear systems of PDEs.
\section{Wavelet Theory}
\label{sec:wavelet}
In this section, we briefly summarize the wavelet theory for clarity of presentation. More information on the wavelet formulation can be found in \cite{goedecker2009wavelets, donoho1992interpolating, beylkin1997adaptive, harnish2023adaptive}.
\subsection{Multiresolution Analysis}
A multiresolution analysis (MRA), which provides a mathematical foundation for the wavelet theory, describes sequential approximation spaces and their duals, denoted by $\mathbf{V}_j$ and $\tilde{\mathbf{V}}_j$, respectively \cite{daubechies1992ten}. The wavelet spaces and their duals, $\mathbf{W}_j$ and $\tilde{\mathbf{W}}_j$, are complements of the approximation spaces in $\mathbf{V}_{j+1}$. These spaces are subsets of the corresponding spaces at higher resolution levels and have the properties
\begin{align}
    \mathbf{V}_j & \subset \mathbf{V}_{j+1}, & \tilde{\mathbf{V}}_j &\subset \tilde{\mathbf{V}}_{j+1}, & \mathbf{V}_{\infty} &= L^2(\Omega), & \mathbf{V}_{j+1} = \mathbf{V}_{j} \oplus \mathbf{W}_j,
\end{align}
where $\Omega$ is some finite domain and $j\in \mathbb{Z}$ denotes resolution level. These nested spaces provide the mathematical foundation for wavelet families of basis functions and their duals. 
The so-called scaling functions ${}^0\phi^j_k$  along with their duals ${}^0\tilde{\phi}^j_k$ are the basis functions in $\mathbf{V}$ and $\tilde{\mathbf{V}}$ spaces respectively. Analogously,  ${}^0\psi^j_k$  and ${}^0\tilde{\psi}^j_k$ are the wavelet basis functions in $\mathbf{W}$ and $\tilde{\mathbf{W}}$ spaces. The index $k\in \mathbb{Z}$ is a local index that provides information that partially defines a location in $\Omega$. These spaces and basis functions have the useful property that their multi-dimensional representations are the tensor product of the one-dimensional components,
\begin{align}
    {}^0\Psi^j_{\vec{k}}(x_1,x_2) &= {}^0\phi^j_{k_1}(x_1) \otimes {}^0\phi^j_{k_2}(x_2).
\end{align}
Wavelet basis functions can be refined and expressed at a higher resolution level as
 \begin{align}
     \phi^{j+1}_{2k}(x) = \phi_k^j(2x - 2k\Delta x^{j+1}),
 \end{align}
 where $\Delta x^{j+1}$ is the grid spacing at level $j+1$ in the direction of $x$. Although all wavelet-based approximations share a common form, there are many unique families of wavelets with different properties \cite{karoui1994families}. The biorthogonal wavelet family used in this work is known as the Deslauriers-Dubuc family \cite{de2003dubuc, goedecker2009wavelets, harnish2018adaptive, harnish2021multiresolution}, with boundary-modified wavelets \cite{vasilyev2000second} used along the spatial and temporal boundaries of the domain.

 \subsection{Spacetime Wavelet Discretization}
We define a finite spacetime domain $\Omega  = \Omega_x \ \times \ [0, T] $ where $\Omega_x \subset \mathbb{R}^N$ and $t \in [0, T]$, $[0, T] \subset \mathbb{R}^{+}$ with spatial boundary $\partial \Omega_x \subset \mathbb{R}^N$. For the problems in this work, we use a 1D spatial domain $\Omega_x$ ($n=1$), which yields a 2D spacetime problem. We generalize all points $k_x$, $k_t$ in the spacetime domain into an array $\vec{k}$ with size governed by the spatial and temporal basis orders, $p_x$ and $p_t$, respectively. Next, we define a vector $\vec{q}$ containing both spatial and temporal variables to simplify notation as $\vec{q} = [x, t]$. The 2D spacetime wavelet representation of a function is expressed as
\begin{align}
    f^j(x, t) = f^j\big(\vec{q}\big) = \sum_{k_i \in [0, 2p_i]} {}^0\mathbb{d}^1_{\vec{k}} \ {}^0\Psi^1_{\vec{k}}(\vec{q}) + \sum_{j=1}^{j_{\mathrm{max}}} \sum_{\lambda = 1}^{2^{N+1}-1} \sum_{\vec{k}}{}^{\lambda}\mathbb{d}^j_{\vec{k}} \ {}^{\lambda}\Psi^j_{\vec{k}}(\vec{q}),
    \label{eq:2dFullapprox}
\end{align}
where index $\lambda \in \mathbb{Z}$ describes the type and location of a wavelet coefficient and $j_{\mathrm{max}}$ is the highest resolution level used. When solving problems on a sparse wavelet grid, the second summation term is necessary because some ${}^{\lambda}\mathbb{d}$ wavelet coefficients can be discarded if determined to be below a user-defined threshold \cite{harnish2023adaptive}. In this work, we use a dense scaling function representation \cite{goedecker2009wavelets} with an expression equivalent to Eq. (\ref{eq:2dFullapprox})
\begin{align}
    f^j\big(\vec{q}\big) = \sum_{k_i \in [0, 2^jp_i]} {}^0\mathbb{d}^j_{\vec{k}} \ {}^0\Psi^j_{\vec{k}}\big(\vec{q}\big).
    \label{eq:2dapprox}
\end{align}
Using this expression simplifies the approximation in that only one type of coefficient, ${}^0\mathbb{d}$, is required to represent the function.
The ${}^0\mathbb{d}$ coefficients can be computed exactly by evaluating the integral with the dual basis
\begin{align}
    {}^0\mathbb{d}^j_{\vec{k}} &= \int_{\Omega} f\big(\vec{q}\big)\tilde{\Psi}^j_{\vec{k}}\big(\vec{q}\big)\mathrm{d}\Omega,
    \label{eq:dCoef}
\end{align}
due to properties of the Deslauries-Dubuc wavelet family \cite{de2003dubuc}. 

The forward wavelet transform (FWT) maps field value coefficients to their corresponding wavelet coefficients and is referred to as wavelet analysis. The inverse process is known as the backward wavelet transform (BWT) or wavelet synthesis, where wavelet coefficients are mapped back to their respective field value coefficients. The operators used to perform the forward and backward transforms are denoted by $\mathbbb{F}$ and $\mathbbb{B}$ where $\mathbbb{B} = \mathbbb{F}^{-1}$.
For the context of this work, it is useful to express the discrete forward and backward wavelet transforms in two dimensions using indicial notation
\begin{align}
    {}^\lambda\mathbb{d}_{mn} &= \mathbb{F}_{mo} \mathcal{F}_{or} \mathbb{F}_{rn}, & \mathcal{F}_{mn} &= \mathbf{\mathbb{B}}_{mo} {}^\lambda\mathbb{d}_{or} \mathbb{B}_{r n},
    \label{eq:analSyn}
\end{align}
where $\mathcal{F}$ and ${}^\lambda\mathbbb{d}$ are the arrays of values in physical and wavelet spaces, respectively. $\mathcal{F}_{mn}$ is populated by evaluating the function $f(\vec{q})$ at each point on the dense spacetime grid. The index $m$ describes the spatial coordinate and $n$ describes the specific time. It is important to note that when solving in the physical space, the coefficients ${}^0\mathbb{d}$ in Eq. (\ref{eq:dCoef}) are equivalent to the field values of $f$ due to the properties of the scaling basis functions \cite{harnish2018adaptive, harnish2021multiresolution, goedecker2009wavelets}. Therefore, we drop the superscript $0$ for brevity of presentation.

\subsection{2D Wavelet Derivative}
In this section, we discuss the wavelet discretization of the $\alpha$th-order derivative of a function taken with respect to a general spacetime variable $q_i$. Applying the $\frac{\partial^{(\alpha)}}{\partial q_i^{(\alpha)}}$ operator to Eq. (\ref{eq:2dapprox}) yields 
\begin{align}
     \frac{\partial^{(\alpha)}f\big(\vec{q}\big)}{\partial q_i^{(\alpha)}}\approx \sum_{\vec{k}} \mathbb{d}^j_{\vec{k}} \frac{\partial^{(\alpha)}\Psi^j_{\vec{k}}\big(\vec{q}\big)}{\partial q_i^{(\alpha)}}.
\end{align}
Because $\Psi\big(\vec{q}\big)$ in the Deslaurier-Dubuc family is continuous and differentiable, its derivative can be represented as a wavelet approximation with new coeffiecients $\Gamma$. These coefficients are commonly referred to as connection coefficients \cite{beylkin1997adaptive}, which again are computed by integrating with the dual basis functions
\begin{align}
    \frac{\partial^{(\alpha)}\Psi_{\vec{k}}^j\big(\vec{q}\big)}{\partial q_i^{(\alpha)}} & =\sum_{\vec{l}} \Gamma^{\vec{l}, j}_{\vec{k}}\Psi_{\vec{l}}^j\big(\vec{q}\big), & \Gamma^{\vec{l}, j}_{\vec{k}} & =\int\left[\frac{\mathrm{d}}{\mathrm{d}q_i}\Psi_{\vec{k}}^j\big(\vec{q}\big)\right]\tilde{\Psi}_{\vec{l}}^j\big(\vec{q}\big)\mathrm{d}q_i,
\end{align} 
where $l \in \mathbb{Z}$ partially describes the location of the derivative coefficients on the computational grid.
The process of computing the connection coefficients reduces to an eigenvector problem \cite{goedecker2009wavelets}. It should be noted that the interpolation order of the bases (\textit{i.e.}, $p_x$ and $p_t$) must be chosen such that the basis functions are sufficiently differentiable for the highest order derivative present in the PDE. The continuity of the Deslauries-Dubuc wavelet family for relevant $p$ values has been studied in \cite{rioul1992simple} and is presented in Table \ref{tab:cont}.
\begin{table} [H]
\centering
    \begin{TAB}(b, 0.1cm, 0.1cm){|c|c|c|}{|c|c|c|c|}
    $p$ & Hölder Regularity & Continuity \\
    4 & $\dot{C}^2$ & $C^1$ \\
    6 & $\dot{C}^{2.830}$ & $C^2$ \\
    8 & $\dot{C}^{3.551}$ & $C^3$ \\
    \end{TAB}
    \caption{Regularity and continuity of Deslaurier-Dubuc wavelets as a function of basis order \cite{rioul1992simple}.}
    \label{tab:cont}
\end{table}

The connection coefficients can be computed for any finite domain, resolution level, and derivative order (so long as a sufficiently-high basis order is selected), and with respect to any independent variable defined on a spacetime domain. ${}^{(\alpha, q_i)}\boldsymbol{\Gamma}$ denotes the $\alpha$th-order derivative operator with respect to $q_i$. The discrete approximation of the derivative is
\begin{align}
    \frac{\partial^{\alpha} f^j(\vec{q})}{\partial q_i^{\alpha}} = {}^{(\alpha,q_i)}\Gamma_{mo}^j \mathbb{d}_{on}^j. \label{eq:dd}
\end{align}
The error of the dense wavelet derivative approximation is described as a function of the dense grid spacing $\Delta q_i$ \cite{dubos2013conservative}
\begin{align}
    \norm{ \frac{\partial^{\alpha} f\big(\vec{q}\big)}{\partial q_i^{\alpha}} - \frac{\partial^{\alpha} f^j\big(\vec{q}\big)}{\partial q_i^{\alpha}}}_{\infty} \leq \mathcal{O}\left(\Delta q_i^{p-\alpha}\right).
    \label{eq:derErr}
\end{align}
As we will show in Section \ref{sec:numerical}, when solving PDEs, the expected order of convergence is defined by the smallest $p-\alpha$ value present in the direction of the largest $\Delta q_i$ on the grid. For example, if the largest grid spacing is in the $q_i$ direction, a second-order PDE ($\alpha = 2$) in $q_i$ solved using $p_{q_i} = 6$ should achieve $4$th-order convergence on a dense wavelet grid. Similar \textit{a priori} error estimates for the function value and nonlinear products have been discussed in previous works \cite{harnish2018adaptive, harnish2023adaptive}.
\section{Computational Implementation}
\label{sec:nonlinearComp}
Computational implementation is performed using the Multiresolution Wavelet Toolkit (MRWT) \cite{harnish2018adaptive, harnish2021multiresolution, harnish2023adaptive}. Due to the nonlinear nature of the problems in this work, we use the Newton-Raphson method to linearize the governing equations and obtain a system of algebraic equations. With the Newton-Raphson method, it is crucial to intelligently select an initial guess to begin the iterations, otherwise convergence may be unattainable. To ensure our method begins within the region of convergence, we propose a wavelet-based recursive solution technique. This recursive strategy also improves the system sensitivity. We solve the algebraic system of equations using the MUltifrontal Massively Parallel sparse direct Solver (MUMPS) \cite{amestoy2000mumps}, which utilizes parallel programming to achieve a faster time to solution. As the problem size increases, so does the computational cost of building and solving the system. Therefore, we also parallelize tangent matrix assembly. 

\subsection{Burgers' Equation}
To illustrate the discretization and implementation of a time-dependent nonlinear PDE using the spacetime wavelet method, we begin by discussing the 1D Generalized Burgers' equation \cite{tersenov2010generalized}. We consider its general form expressed as
 \begin{align}
    \begin{split}
         &\frac{\partial f(\vec{q})}{\partial t} + \biggl(f(\vec{q}) +c(t)\biggl)\frac{\partial f(\vec{q})}{\partial x} - \nu \frac{\partial ^2 f(\vec{q})}{\partial x^2}  = 0 \ \ \mathrm{in} \ \ \Omega, \\
         &f(\vec{q}) = f_I  \ \ \mathrm{on} \ \ \Omega_x \times (t = 0), \\ 
         &f(\vec{q}) = f_B  \ \ \mathrm{on} \ \ \partial\Omega_x \times (0, T],
     \end{split}
     \label{eq:modBurg}
 \end{align}
where $c(t)$ and $\nu$ are coefficients of advection and diffusion, respectively. The discrete wavelet residual for the Burgers' equation is expressed in matrix form as
\begin{align}
    \mathbf{R} \ = \ \boldsymbol{\hat{\mathcal{F}}} \cdot {}^{(1, t)}\boldsymbol{\Gamma} \ + \ 
    \left(\boldsymbol{\hat{\mathcal{F}}} +\mathbf{c} \right)\circ \left({}^{(1, x)}\boldsymbol{\Gamma} \cdot \boldsymbol{\hat{\mathcal{F}}} \right) \ - \ \nu \ {}^{(2, x)}\boldsymbol{\Gamma}\cdot\boldsymbol{\hat{\mathcal{F}}} = \mathbf{0},
    \label{eq:mb0}
\end{align} 
where $\boldsymbol{\hat{\mathcal{F}}} = \boldsymbol{\mathcal{F}_U} + \boldsymbol{\mathcal{F}}_B + \boldsymbol{\mathcal{F}}_I$ is an array containing all degrees of freedom, including boundary conditions $f_B$, and initial condition $f_I$. $\boldsymbol{\mathcal{F}}_U$ is the array of unknowns. The $\circ$ symbol denotes a Hadamard product. The derivative terms are derived from Eq. (\ref{eq:dd}). Note that the temporal derivative operator ${}^{(\alpha, t)}\boldsymbol{\Gamma}$ contracts with the second index of $\boldsymbol{\hat{\mathcal{F}}}$, which corresponds to the time, and ${}^{(\alpha, x)}\boldsymbol{\Gamma}$ contracts with the spatial index.
\subsection{Newton-Raphson Method} \label{sec:NR}
Applying Newton's method to the Burgers' equation, we use the discretized residual given by Eq. (\ref{eq:mb0}) and linearize it to obtain the Jacobian
 \begin{align}
\frac{\partial \mathbf{R}}{\partial \boldsymbol{\hat{\mathcal{F}}}} = \mathbf{K} = \mathbf{I}\otimes{}^{(1, t)}\boldsymbol{\Gamma} \ + \ 
\left({}^{(1, x)}\boldsymbol{\Gamma}\cdot\boldsymbol{\hat{\mathcal{F}}} \right)\otimes\mathbf{I} \ + \ 
\left(\boldsymbol{\hat{\mathcal{F}}} + \mathbf{c}\right) \circ {}^{(1, x)}\boldsymbol{\Gamma}\otimes\mathbf{I} \ - \ 
\nu \ {}^{(2, x)}\boldsymbol{\Gamma}\otimes\mathbf{I},
\label{eq:tan}
\end{align}
where $\mathbf{I}$ denotes the identity matrix.

We vectorize the corrections $\boldsymbol{\Delta \mathcal{F}}$ and residuals by stacking their columns, and we partition the tangent matrix $\mathbf{K}$ to separate the unknown interior entries $\mathbf{K}_{UU}$ from the known boundary and initial conditions:
\begin{align}
\begin{bmatrix}
\mathbf{K}_{UU} & \mathbf{K}_{UB} & \mathbf{K}_{UI}\\
\mathbf{K}_{BU} & \mathbf{K}_{BB} & \mathbf{K}_{BI} \\
\mathbf{K}_{IU} & \mathbf{K}_{IB} & \mathbf{K}_{II}
\end{bmatrix}
\begin{bmatrix}
\Delta\vec{\mathcal{F}}_U \\ 
\Delta\vec{\mathcal{F}}_B \\ 
\Delta\vec{\mathcal{F}}_I
\end{bmatrix}
\ = \ -
\begin{bmatrix}
\vec{\mathcal{R}}_U \\ 
\vec{\mathcal{R}}_B \\ 
\vec{\mathcal{R}}_I
\end{bmatrix}.
\label{eq:kMat}
\end{align}
This technique yields the system of algebraic equations, 
\begin{align}
    \mathbf{K}_{UU}\Delta\vec{\mathcal{F}}_U = -\left(\vec{\mathcal{R}}_U + \mathbf{K}_{UI}\Delta\vec{\mathcal{F}}_I + \mathbf{K}_{UB}\Delta\vec{\mathcal{F}}_B\right) = -\vec{\mathcal{R}}
    \label{eq:bcSys}
\end{align}
to be solved for unknown corrections $\Delta\vec{\mathcal{F}}_U$. As is typical, we update our solution as $\vec{\mathcal{F}}^{i+1}_U = \vec{\mathcal{F}}^{i=0}_U + \Delta\vec{\mathcal{F}}^{i}_U$, where $\Delta \vec{\mathcal{F}}^i_U = \sum_{i} \Delta \Delta \vec{\mathcal{F}}^i_U$. These corrections are made until a user defined tolerance, $\sfrac{\norm{\mathcal{R}}}{\norm{\mathcal{R}_{0}}} < \mathrm{tol}$, is obtained. We have now enforced boundary and initial conditions and are left with a system of algebraic equations resembling that of an elliptic and/or implicitly discretized time-dependent PDE. We note that the initial and boundary terms in Eq. (\ref{eq:bcSys}) are applied only in the first iteration $i=0$ and the subsequent $\Delta\Delta \vec{\mathcal{F}}_I$ and $\Delta\Delta \vec{\mathcal{F}}_B$ corrections are equal to zero \cite{simo2006computational}. One can see that we have now introduced two forcing terms: $\mathbf{K}_{UI}\Delta\vec{\mathcal{F}}_I $ resulting from the initial condition, and $\mathbf{K}_{UB}\Delta\vec{\mathcal{F}}_B$ from the boundary conditions. Depending on the nature of the initial and boundary conditions, these forcing terms can vary significantly in magnitude. Additionally, the conditioning of the tangent matrix $\mathbf{K}_{UU}^{i=0}$ can be dependent on the selected initial guess $\vec{\mathcal{F}}^{i=0}_U$ and the resolution level $j$.

\subsection{Wavelet-Based Recursive Solution Technique}
\label{sec:recursive}
Due to the implicit nature of our systems, a poorly selected initial guess and/or poor matrix conditioning will lead to a divergent solution. Additionally, we will show that the magnitude of the right-hand side $\lVert\vec{\mathcal{R}}\rVert$ also impacts the solution sensitivity. Consider our system in Eq. (\ref{eq:bcSys}) and apply a perturbation $\delta \vec{\mathcal{R}}$ to the right-hand side, resulting in a perturbation to the solution $\delta \vec{\mathcal{F}}$, yielding the system
\begin{align}
    \mathbf{K}_{UU}\left(\Delta \vec{\mathcal{F}}_U + \delta \vec{\mathcal{F}}\right) = -(\vec{\mathcal{R}} + \delta \vec{\mathcal{R}}).
\end{align}
Subtracting the original system, applying the matrix inverse, and taking norms, we obtain 
\begin{align}
    \lVert\delta \vec{\mathcal{F}}\rVert = \lVert\mathbf{K}_{UU}^{-1} \delta \vec{\mathcal{R}}\rVert \leq \lVert\mathbf{K}_{UU}^{-1}\rVert\lVert\delta \vec{\mathcal{R}}\rVert.
    \label{eq:normP}
\end{align}
Recombining Eq. (\ref{eq:normP}) with the norm of the original problem, we manipulate the equation to obtain the typical expression describing the sensitivity of an algebraic system \cite{heath2018scientific}
\begin{align}
    \frac{\lVert\delta \vec{\mathcal{F}}\rVert}{\lVert\vec{\Delta \mathcal{F}}_U\rVert} \leq \kappa \left(\mathbf{K}_{UU}\right)  \frac{\lVert\delta \vec{\mathcal{R}}\rVert}{\lVert\vec{\mathcal{R}}\rVert},
    \label{eq:solvability}
\end{align}
where $\kappa \left(\mathbf{K}_{UU}\right)$ denotes the condition number of the tangent matrix. From this expression, we see that two factors impact the sensititivity of the solution: the condition number of the matrix $\mathbf{K}_{UU}$ and the magnitude of the right-hand side $\vec{\mathcal{R}}$. 

Unlike traditional time-marching explicit or implicit methods \cite{mosbach2009quantitative, verwer1996explicit, moin2010fundamentals}, the spacetime wavelet solver is not subject to a stability condition. However, motivated by Eq. (\ref{eq:solvability}), we define a \textit{solvability condition}, which is a function of wavelet basis order, initial guess, $\kappa\left(\mathbf{K}_{UU}\right)$ and $\lVert \vec{\mathcal{R}}\rVert$. We define a system as solvable if $\Delta \Delta \vec{\mathcal{F}}^i_U \rightarrow 0$ as $i \rightarrow \infty$. We note that the basis order is selected to accommodate the highest derivative in the system and is fixed during the analysis.
Therefore, we focus on controlling the initial guess. A better initial guess leads to a smaller residual vector $\lVert\vec{\mathcal{R}}\rVert$ making the system less sensitive according to Eq. (19). Moreover, it maintains favorable tangent conditioning. Using Eq. (\ref{eq:solvability}) as a guide, we develop a wavelet-based iterative procedure to reduce the sensitivity of the system by simultaneously adjusting $\kappa\left(\mathbf{K}_{UU}\right)$ and $\lVert \vec{\mathcal{R}}\rVert$.


Here we describe the procedure to improve solvability of the system by constructing a wavelet-guided  initial guess.  Moreover, we utilize an \textit{a priori} error estimate in Eq. (\ref{eq:derErr}) to control the accuracy of the spacetime solution.  Traditional time-marching methods select the timestep  $\Delta t$  as $\Delta t = \mathrm{min} (\Delta t_s , \Delta t_a )$, where $\Delta t_s$ is the minimum timestep required to obtain a stable solution and $\Delta t_a$ is the minimum timestep required to obtain a well-resolved (\textit{i.e.}, accurate) PDE solution.  In general, time-marching schemes control accuracy of the solution \textit{a posteriori}. In our work, we utilize the \textit{a priori} error estimate in Eq. (\ref{eq:derErr}) that provides the spacing $\Delta q_i$ at a resolution level $j_\mathrm{max}$. Then the accuracy of the solution is given by the largest wavelet coefficient on level $j_{\mathrm{max}}+1$ \cite{harnish2023adaptive}. However, when solving a large system of algebraic equations at $j_{\mathrm{max}}$ with an uninformed initial guess (\textit{e.g}., $\vec{\mathcal{F}}^{i=0}_U = \vec{0}$), the problem is often outside of the region of convergence and/or can be poorly conditioned.
Therefore to obtain a convergent Newton-Raphson method and satisfy the accuracy requirements given by Eq (\ref{eq:derErr}), we use MRA to build the solution recursively. 

\begin{figure}[H]
    \centering
    \includegraphics[width=1\textwidth]{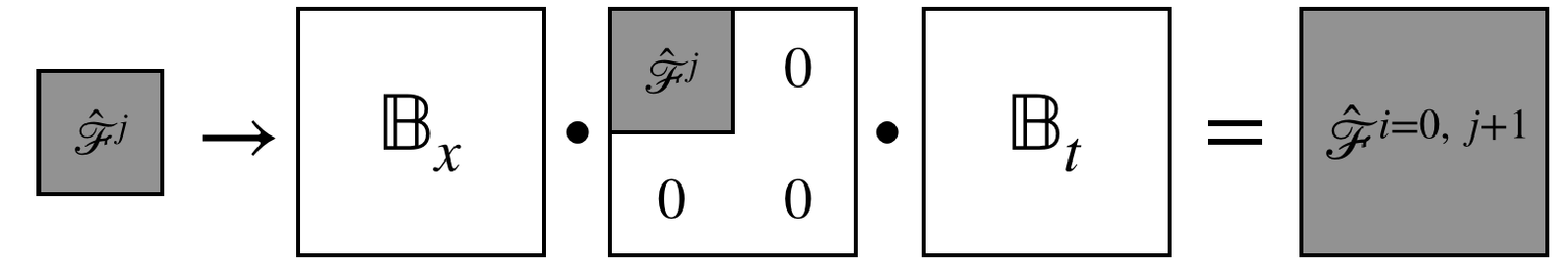}
    \caption{Wavelet synthesis of initial guess $\hat{\mathcal{F}}^{i = 0, \ j+1}$ from previous solution $\hat{\mathcal{F}}^{j}$.}
    \label{fig:ig}
\end{figure}

We begin at some level $j<j_{\mathrm{max}}$ using a generic initial guess (\textit{e.g.}, zeros). After the first solution is found, we use wavelet synthesis to construct the solution at level $j+1$ to serve as the initial guess for the next iteration and also to check the accuracy of the solution. Moreover, this process impacts the conditioning of $\mathbf{K}_{UU}$ and moves the solution closer to the region of convergence. The process of prolongating the solution with wavelet synthesis is illustrated by Fig. \ref{fig:ig}. We repeat this process until $j_{\mathrm{max}}$ is reached and the desired accuracy given by Eq. (\ref{eq:derErr}) is obtained. Unfortunately, even with a well-conditioned tangent matrix, the $\Delta \Delta \vec{\mathcal{F}}_U^i$ corrections can begin to diverge if $\lVert\vec{\mathcal{R}}\rVert$ is too large and sensitive to perturbations, as illustrated by Eq. (\ref{eq:solvability}) (\textit{i.e.}, solvability condition), forcing iterations to begin outside of the region of convergence. To account for this issue, we simultaneously control the magnitude of the residual vector $\lVert\vec{\mathcal{R}}\rVert$ while improving the conditioning of the system.

The magnitude of the right-hand side of Eq. (\ref{eq:bcSys}) is a function of both the initial guess (affects $\vec{\mathcal{R}}_U$, $\mathbf{K}_{UI}$, and $\mathbf{K}_{UB}$) and the initial and boundary conditions (affects $\Delta \vec{\mathcal{F}_I}$ and $\Delta \vec{\mathcal{F}_B}$). In particular, for problems with sharp gradients such as shocks, the initial and boundary conditions can be problematic. For such cases, we specify a series of modified initial and/or boundary conditions for all levels  $j$ to $j_{\mathrm{max}}$. It is important that this series of initial/boundary conditions sufficiently reduces the sensitivity of the system, Eq. (\ref{eq:bcSys}), by changing the magnitudes of $\Delta \vec{\mathcal{F}}_I$ and $\Delta \vec{\mathcal{F}}_B$ in such a way that balance is achieved with $\kappa \left(\mathbf{K}_{UU}\right)$. In this way we satisfy the solvability condition. If we wish to perform $M$ solves to recursively reduce the magnitude of $\lVert\vec{\mathcal{R}}\rVert$ on the system (see Eq. (\ref{eq:solvability}), we define length-$M$ vectors of weights $\vec{\delta} = [\delta_1, \ \delta_2, \ \dots , \ \delta_{M}]$ and  $\vec{\chi} = [\chi_1, \ \chi_2, \ \dots , \ \chi_{M}]$, which specify series of modified initial and boundary conditions $f_{I, g} = f_I(\vec{x}, t_0, \delta_g)$, $f_{B, g} = f_B(\partial \Omega_x, t, \chi_g)$, respectively. These modified conditions recursively grow closer in value to the actual initial and boundary conditions of the problem, satisfying the  solvability condition for each problem in the sequence.

At the final step of the recursive algorithm, we solve the problem at level $j_{\mathrm{max}}$, with desired initial and boundary conditions. During this step, the initial guess synthesized from the $j_{\mathrm{max}}-1$ solution moves the full problem towards the region of convergence to initialize the final Newton-Raphson loop. Finally, we assess if a user-prescribed tolerance is smaller than our \textit{a priori} error estimate  in Eq. (\ref{eq:derErr}). If the user-given tolerance is not yet met, we increase $j_{\mathrm{max}}$ and repeat the process. The recursive algorithm is outlined in Algorithm \ref{alg:recAlg}. The robustness of Algorithm \ref{alg:recAlg} can be improved by employing a globally-convergent Newton-Raphson method \cite{eisenstat1994globally} to ensure the residual does not increase between Newton-Raphson iterations.
\begin{algorithm}[H]
\caption{Recursive Spacetime Wavelet Solver}\label{alg:recAlg}
\begin{algorithmic}
\State Read Input
\State Define starting level $j$ and initial guess $\vec{\mathcal{F}}_U^{i=0}$
\State Define weighting vectors $\vec{\delta}$, $\vec{\chi}$ and set $g = 1$
\While {$j \leq j_{\mathrm{max}}$}
        \State Set conditions: $f_{I,g} = f_I(x, t_0, \delta_g)$,  $f_{B, g} = f_B(\partial \Omega_x, t, \chi_g)$
        \State Compute sparsity pattern of the tangent matrix
        \While{$\left( \text{error} > \text{tolerance}\right)$}
                \State Update tangent $\mathbf{K}_{UU}$ and residual $\vec{\mathcal{R}}_U$ \comment{Eqs. (\ref{eq:mb0}-\ref{eq:tan})} 
                \State Solve for $\Delta \vec{\mathcal{F}}_U$ and update solution $\vec{\mathcal{F}}_U^{i+1} = \vec{\mathcal{F}}_U^i + \Delta \vec{\mathcal{F}}_U^i$ \comment{Eq. (\ref{eq:bcSys})}
        \EndWhile
        \If {$j < j_{\mathrm{max}}$}
            \State Synthesize initial guess from solution \comment{Eq. \ref{eq:analSyn}}
            \State $g = g+1$
        \EndIf
        \State $j = j+1$
\EndWhile
\If{Error requirement is met} \comment{Eq. (\ref{eq:derErr})}
    \State Compute error and output the solution
\Else 
    \State Increase $j_{\mathrm{max}}$ and restart algorithm
\EndIf
\end{algorithmic}
\end{algorithm} 
\section{Numerical Examples}
In this section, we present numerical verification problems for the spacetime wavelet solver. Included are three variations of the nonlinear Burgers' equation: one modeling simple shock advection, one modeling shock advection with time-dependent advection, and one modeling shock evolution. The Burgers' equation variations have analytical solutions allowing rigorous verification studies. The derivatives for the shock advection form of Burgers' equation also have a closed form, which  allows for a full analysis of the numerical convergence of these problems. Finally, we present the Sod shock tube problem, used to test the capability of the spacetime method to solve a system of coupled nonlinear Navier-Stokes equations. All examples are solved with a relative Newton-Raphson tolerance of $10^{-6}$ (see Section \ref{sec:NR}).
\label{sec:numerical}
\subsection{Burgers' Equation} \label{sec:burgers}
Variations of the nonlinear Burgers' equation are used to model both shock advection and shock evolution, governed by Eq. (\ref{eq:modBurg}). The initial-boundary value problem is solved on domain $x \in[-1, 1]$, $t\in [0,\sfrac{1}{2}]$ with a diffusion coefficient $\nu = 0.01$ and transition region centered about $x= -0.5$. The reported errors for all Burgers' examples are calculated by evaluating the highest wavelet coefficient on $j_{\mathrm{max}} +1$ via wavelet synthesis. For the walking Burgers' problem with constant advection (Section \ref{sec:mb}) and the steepening Burgers (Section \ref{sec:steep}) problem we use the basis orders $p_x = 6$, $p_t = 4$. This choice of basis order combination is made based on the lowest order required to differentiate the basis functions in order to minimize the computational work needed to solve the problems. For the walking Burgers' problem with time-dependent advection (Section \ref{sec:mbct}), we use the combination of $p_x=p_t=6$ in order to fully capture the rapidly-changing solution features. For all Burgers' examples, the solvability condition is satisfied using a zero initial guess so the recursive technique is not fully-utilized. Therefore we set $j = j_{\mathrm{max}}$, and the weighing vectors $\vec{\delta}$ and $\vec{\chi}$ contain only the values describing the unmodified initial and boundary conditions of the problem. We will utilize the recursive spacetime algorithm with the Sod problem in Section \ref{sec:sst}.

\subsubsection{Walking Burgers' Equation with Constant Advection (Shock Advection)} \label{sec:mb}
\begin{figure}[H]
\begin{subfigure}{0.49\textwidth}
  \includegraphics[width=\linewidth]{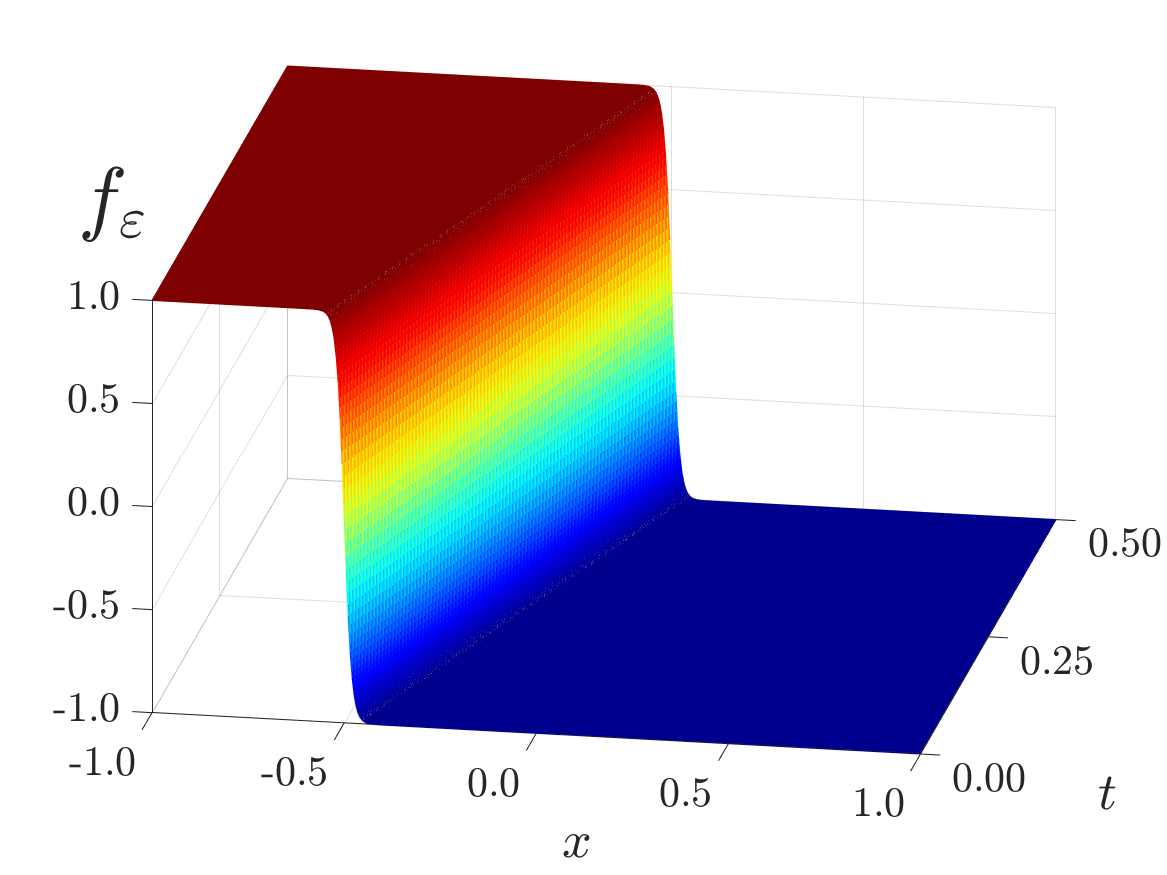}
  \caption{Spacetime solution.}\label{fig:mbST}
\end{subfigure}
\hspace*{\fill}
\begin{subfigure}{0.49\textwidth}%
  \includegraphics[width=\linewidth]{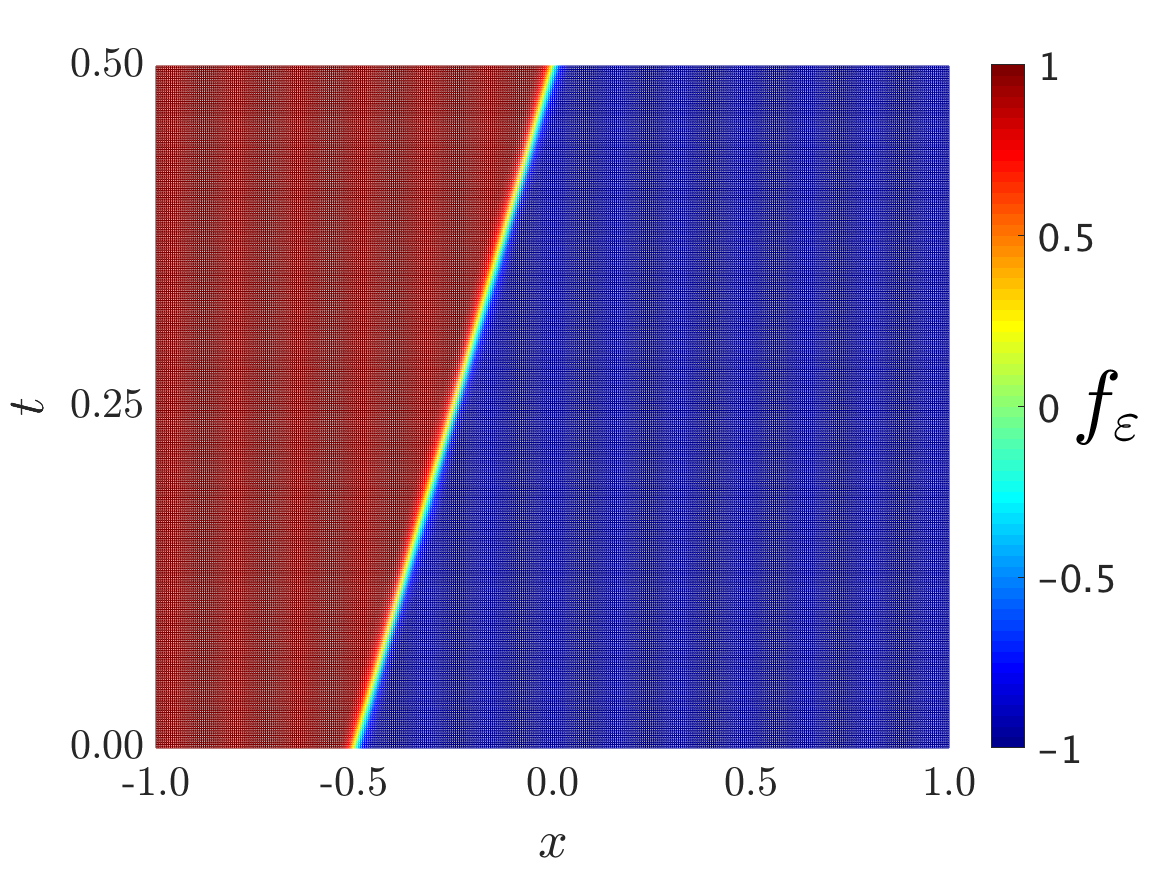}
  \caption{Top view.}\label{fig:mbTop}
\end{subfigure}
\caption{Spacetime solution and profile for shock advection problem with $c(t) = 1$ at $j=7$ \big(Eq. (\ref{eq:modBurg})\big) and $p_x = 6$, $p_t = 4$.}
\label{fig:mb3d}
\end{figure}

\begin{figure}[H]
\begin{subfigure}{0.49\textwidth}
  \includegraphics[width=\linewidth]{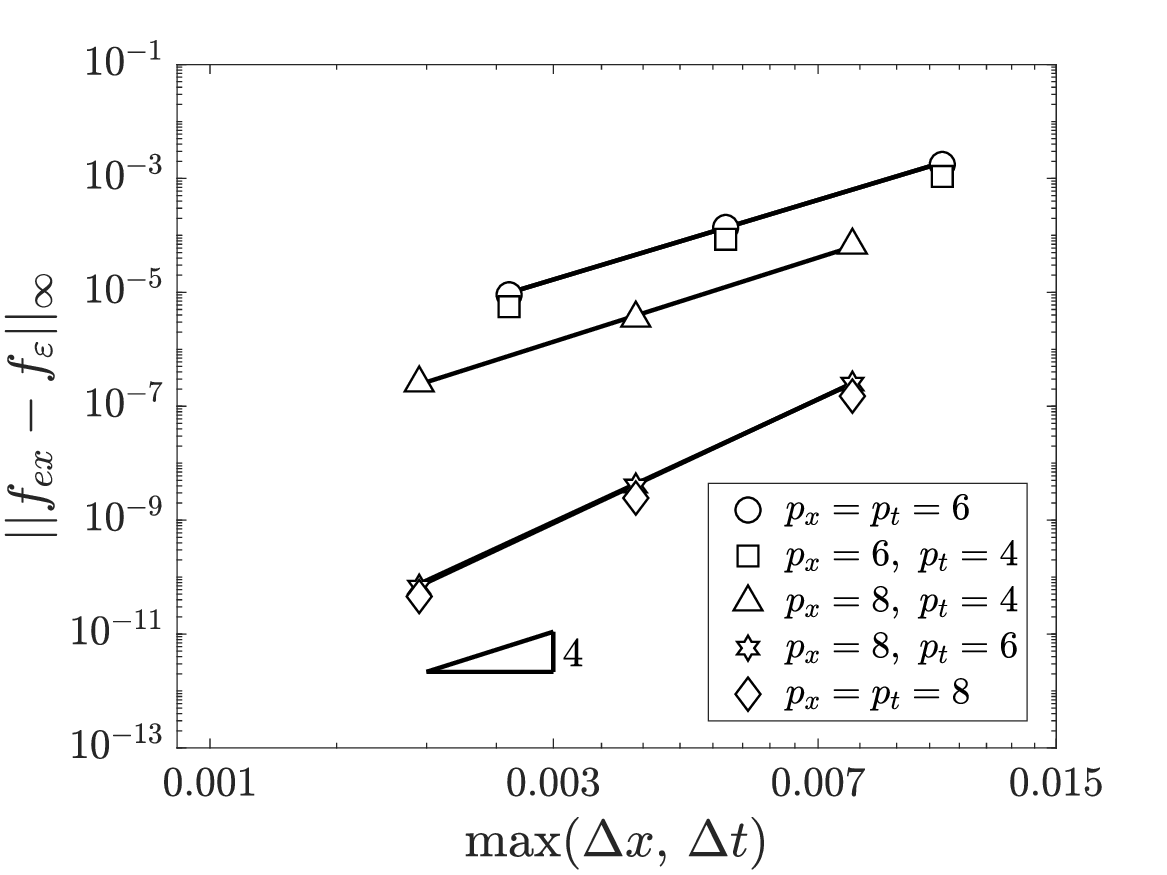}
  \caption{Solution convergence using varying combinations of $p_x$ and $p_t$.}\label{fig:mbSolConv}
\end{subfigure}
\hspace*{\fill}
\begin{subfigure}{0.49\textwidth}%
  \includegraphics[width=\linewidth]{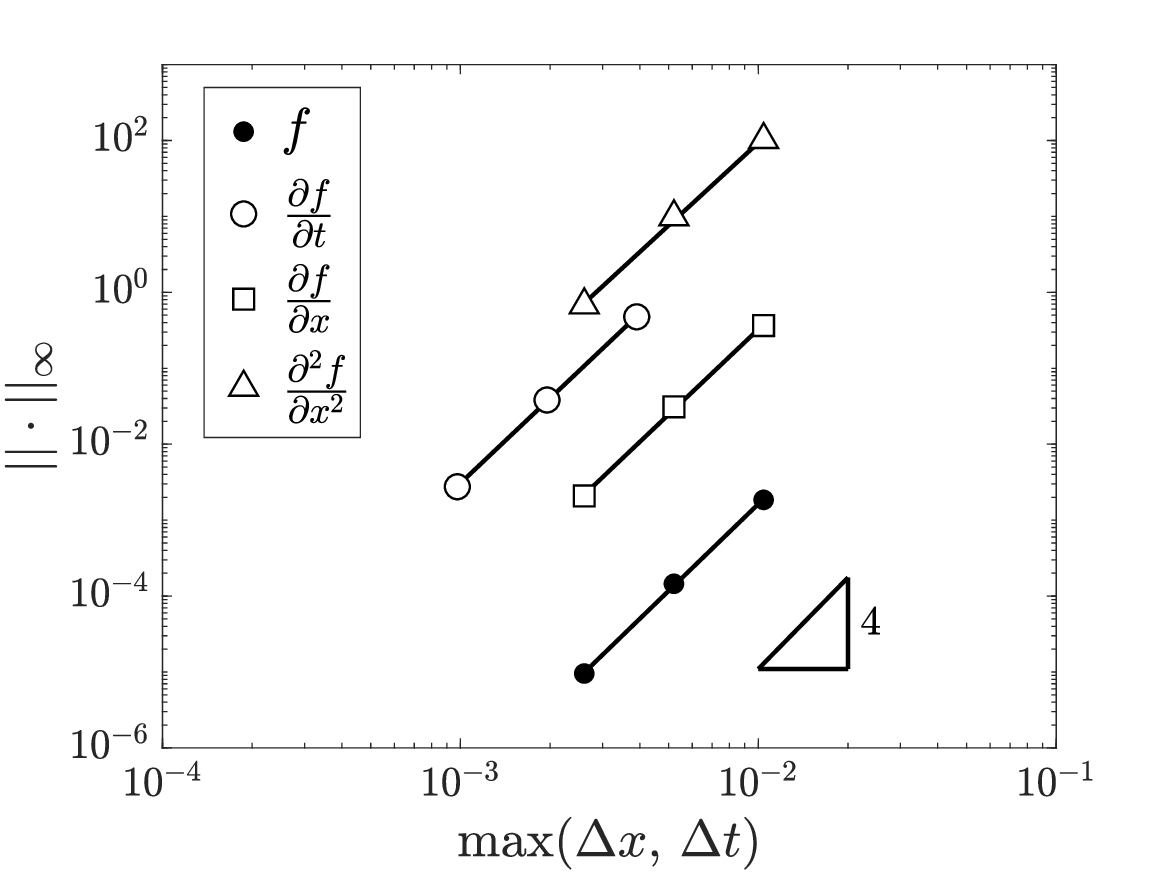}
  \caption{Solution and derivative convergence for $p_x = 6$, $p_t = 4$.}\label{fig:mbDconv}
\end{subfigure}
\caption{Solution and derivative convergence for the shock advection problem with $c(t) = 1$ at $j =5, 6, 7$.}
\label{fig:mbConv}
\end{figure}

The so-called walking Burgers' or Generalized Burgers' equation is defined when the advection term $c(t)$ is nonzero, yielding solution behavior that simulates shock propagation. For the first example, we define $c(t) = 1$. The Dirichlet spatial boundary conditions and initial conditions (\textit{i.e.}, when $t=0$) are obtained from the analytical solution
 \begin{align}
    f(x,t) = -\text{tanh}\left(\frac{x - x_0 - t \ c(t)}{2\nu}\right).
    \label{eq:mbEx}
\end{align}
The resulting spacetime solution for this problem is shown in Fig. \ref{fig:mb3d}. Fig. \ref{fig:mbTop} displays the constant shock advection through the spacetime domain as well as the thickness of the shock governed by the viscosity parameter $\nu$. Fig. \ref{fig:mbConv} shows the spacetime convergence rates for the different combinations of the basis functions.  We can observe that both $p_x=p_t=6$ and $p_x=6$, $p_t=4$ solutions converge with expected $4$th-order rate. The $p_x=8$ and $p_t=4$ solution is also $4$th-order convergent, but the error of the solution is smaller. The $p_x=p_t=8$ and $p_x=8$, $p_t=6$ solutions are the most accurate and converge with higher order (\textit{e.g.}, $p_x=p_t=8$ combinations yield $6$th-order convergence for second-order PDEs).

\begin{figure}[H]
    \centering
    \includegraphics[width=0.5\textwidth]{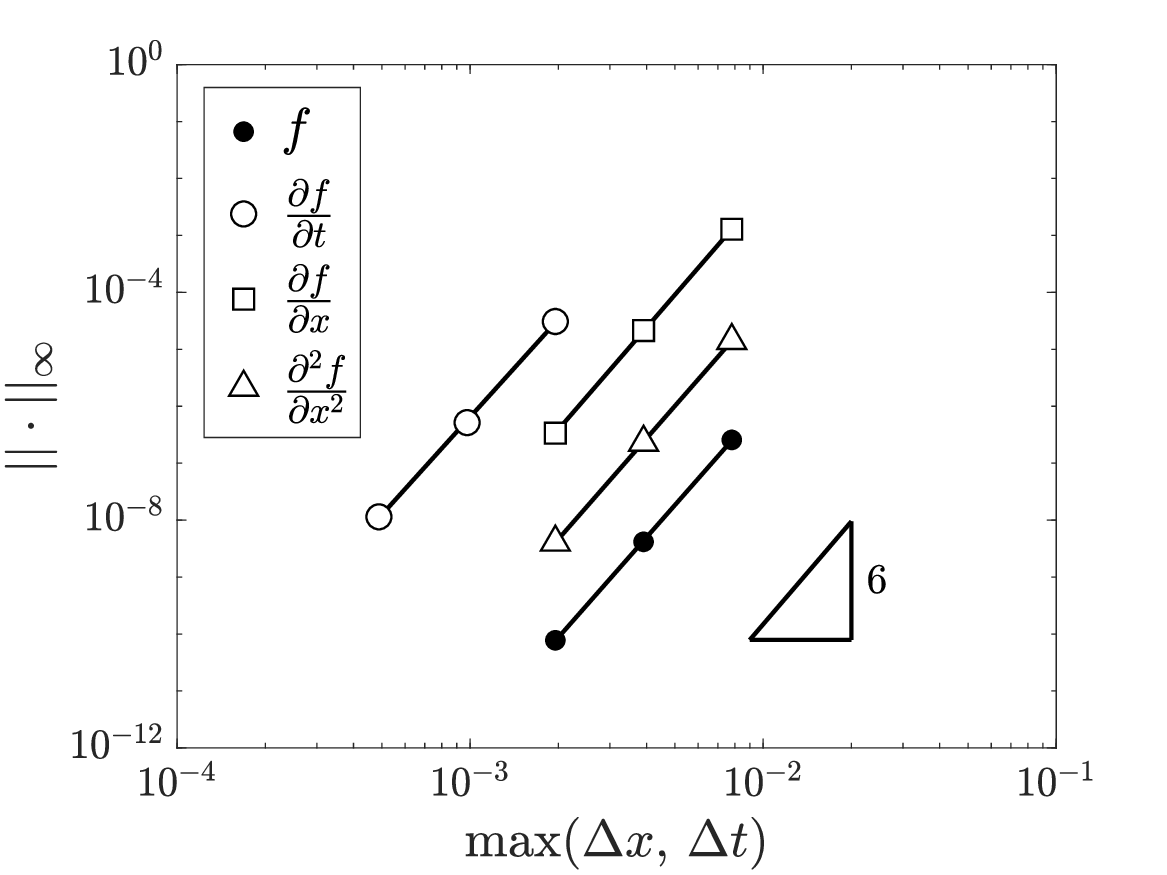}
    \caption{$6$th-order solution and derivative convergence for $j = 5,6,7$ with $p_x=p_t=8$.}
    \label{fig:conv88}
\end{figure}

 A solution accuracy of $\lVert f_{ex} - f_{\varepsilon}\rVert_{\infty} = 9.5 \times 10^{-6}$ is achieved at level $j = 7$ (1,575,425 degrees of freedom (DOF)) along with the expected $4$th-order convergence for both solution and derivative approximations as predicted by Eq. (\ref{eq:derErr}). This makes the scheme fully high-order as we achieve the same order of convergence not only for the field, but also for spatial and temporal derivatives (see Fig. \ref{fig:mbDconv}). As expected, the approximation of the $2$nd-order spatial derivative is the least accurate. We show in Fig. \ref{fig:conv88} that we achieve the expected $6$th-order convergence for both solution and derivative approximation when using basis orders $p_x=p_t=8$, as predicted by Eq. (\ref{eq:derErr}), with a solution accuracy of $\lVert f_{ex} - f_{\varepsilon}\rVert_{\infty} = 7.8 \times 10^{-11}$ at level $j=7$ (4,198,401 DOF).
\noindent These higher-order results serve to demonstrate the capability of the spacetime wavelet solver to match the theoretical convergence rate, Eq. (\ref{eq:derErr}), when increasing the order of the wavelet basis functions.
 
\begin{figure}[H]
    \centering
    \includegraphics[width=0.5\textwidth]{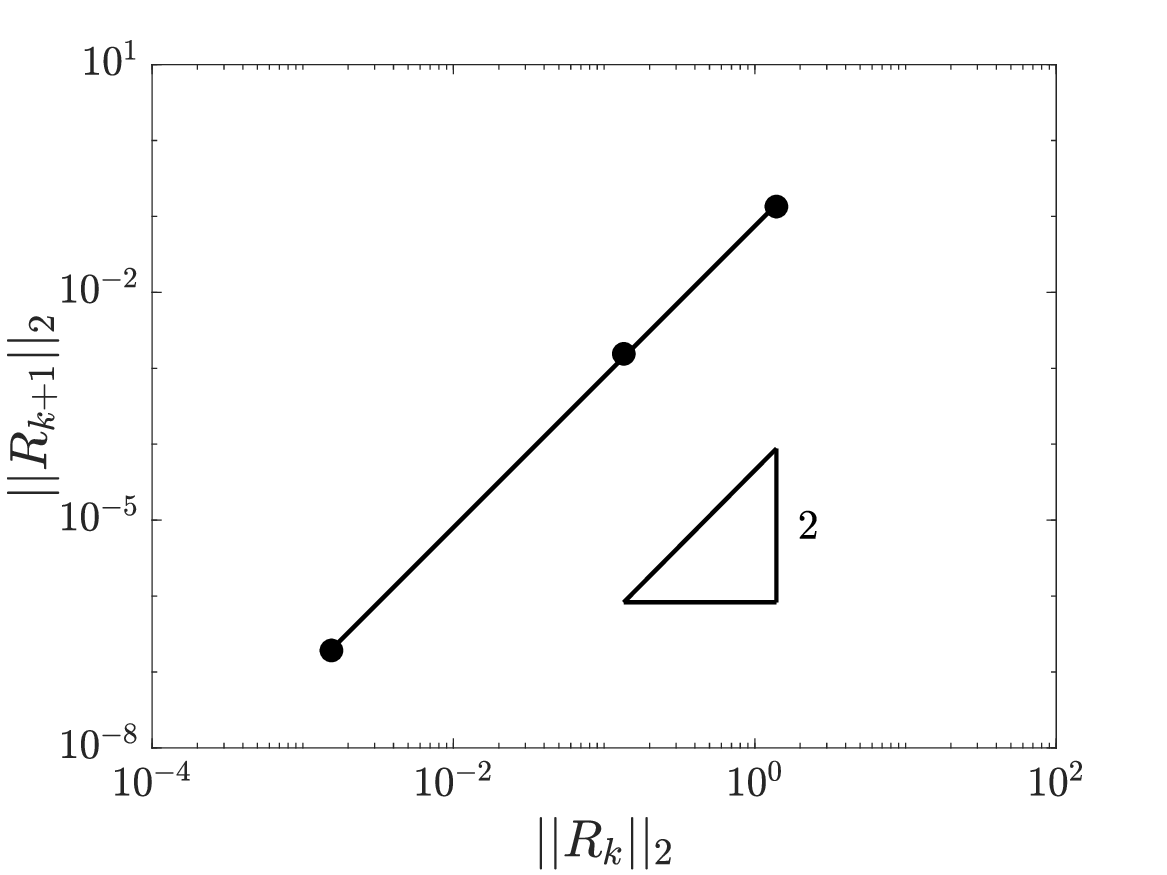}
    \caption{Quadratic convergence rate of Newton-Raphson method while solving the shock advection problem.}
    \label{fig:NR}
\end{figure}

Fig. \ref{fig:NR} displays the quadratic convergence rate of the Newton-Raphson method when applied to the shock advection problem, which indicates the proper consistent linearization of Eq. (\ref{eq:tan}).

\begin{figure}[H]
\begin{subfigure}{0.49\textwidth}
  \includegraphics[width=\linewidth]{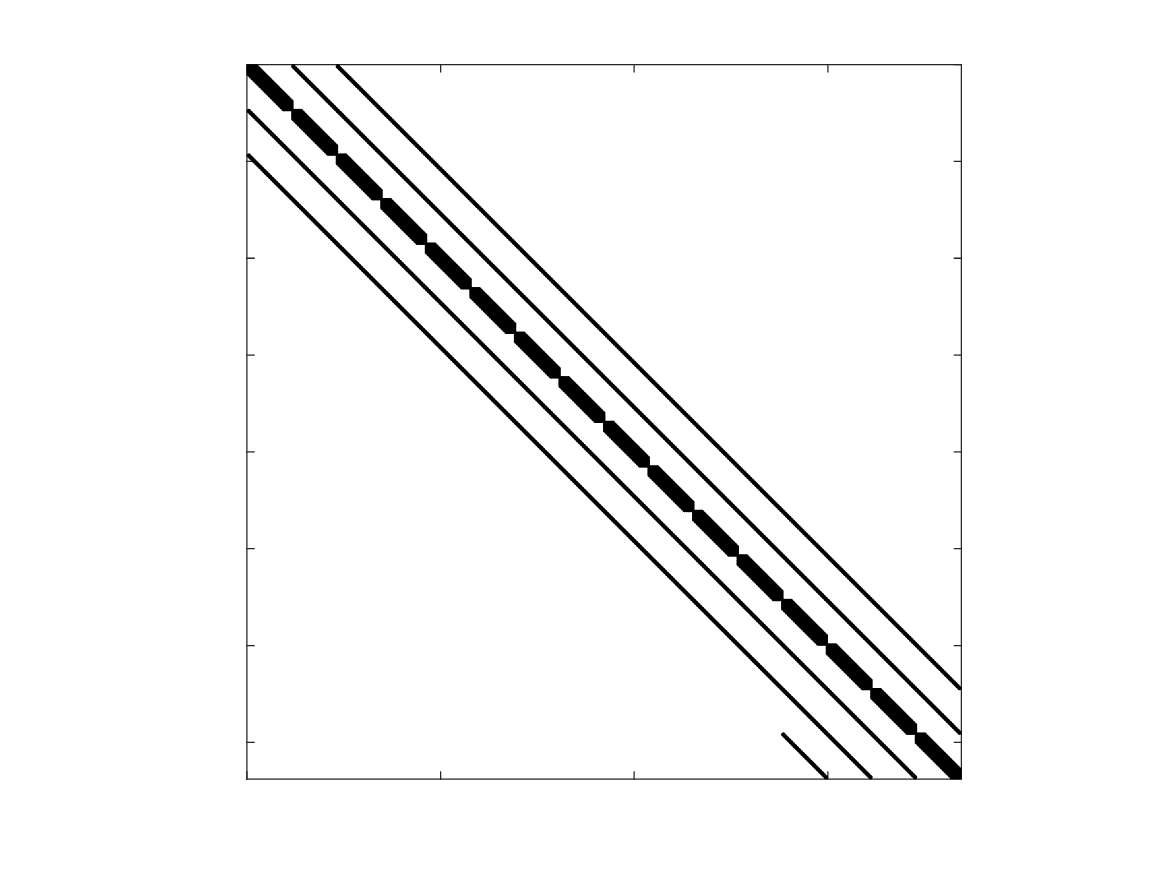}
  \caption{Shock advection sparsity pattern.}\label{fig:mbSparsity}
\end{subfigure}
\hspace*{\fill}
\begin{subfigure}{0.49\textwidth}%
  \includegraphics[width=\linewidth]{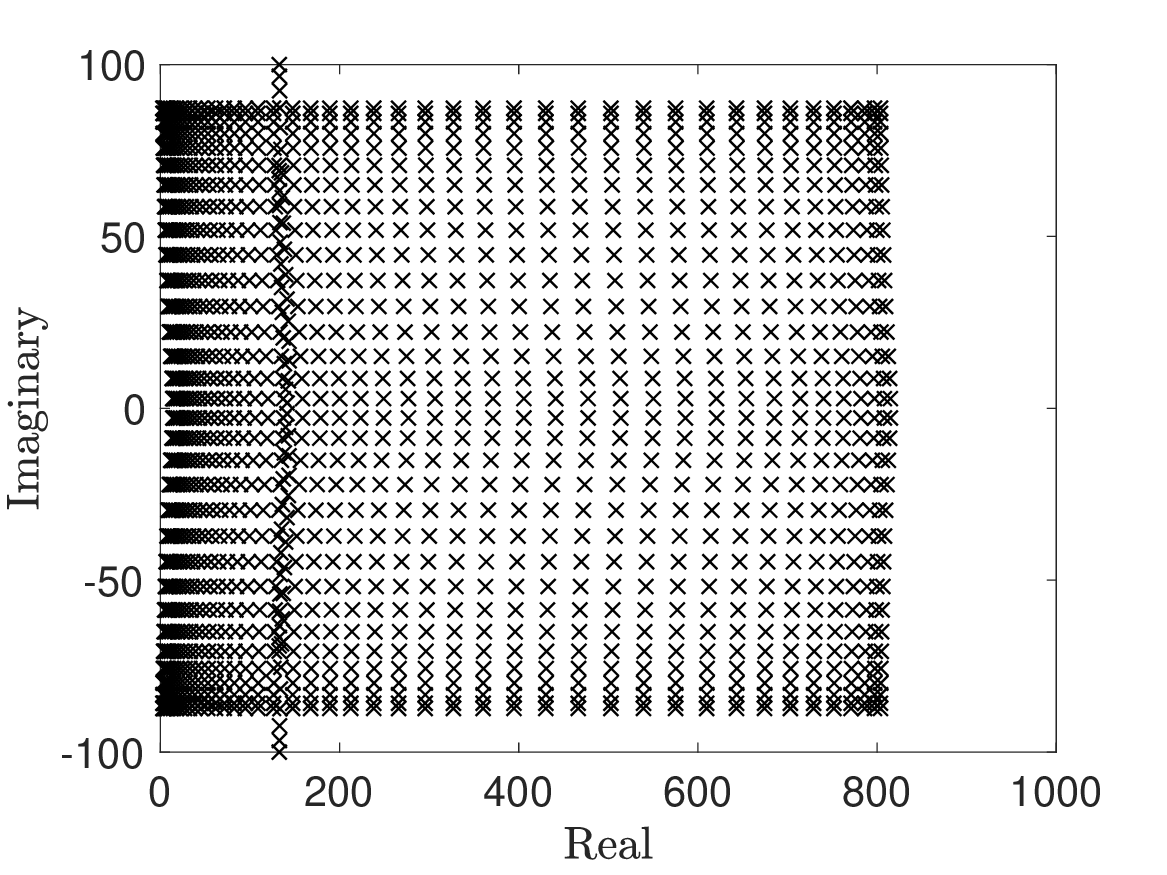}
  \caption{Shock advection eigenvalue spectrum.}\label{fig:mbEigs}
\end{subfigure}
\caption{Spacetime shock advection matrix sparsity pattern and eigenvalue spectrum at $j = 3$ with $p_x = 6$, $p_t=4$.}
\label{fig:mbAnalysis}
\end{figure}

The sparsity pattern and eigenvalue spectrum for the resulting tangent matrix for this problem are shown in Fig. \ref{fig:mbAnalysis}. Note that these plots present the tangent at level $j = 3$ in order to more clearly show relevant features. 
The non-symmetric tangent sparsity pattern has a block diagonal structure of non-zeros with the widest band of non-zeros along the main diagonal. This pattern arises from the Kronecker products of banded derivative operators ($\boldsymbol{\Gamma}$) with the identity matrix present in the Jacobian, Eq. (\ref{eq:tan}). The bandwidth of the main diagonal is determined by the basis order of the first index, $p_x$, and the number of off-diagonal bands is determined by the basis order of the second index, $p_t$. The small band of nonzero entries in the bottom right of Fig. \ref{fig:mbSparsity} results from the linearization of the nonlinear term in Eq. (\ref{eq:mb0}). The eigenvalues (Fig. \ref{fig:mbEigs}) for this problem are all complex with positive real components, characteristic of an elliptic system with smoothed discontinuities \cite{powers2015mathematical}. 

\subsubsection{Comparison to Runge-Kutta Method}
In this section, we compare our spacetime method to the explicit time-marching wavelet solver developed by Harnish et al. \cite{harnish2018adaptive}. Due to its $4$th-order convergence properties and constant timesteps, the $4$th-order Runge-Kutta method (RK4) is a good candidate to compare to the spacetime wavelet method.

\begin{figure}[H]
    \centering
    \includegraphics[width=0.5\textwidth]{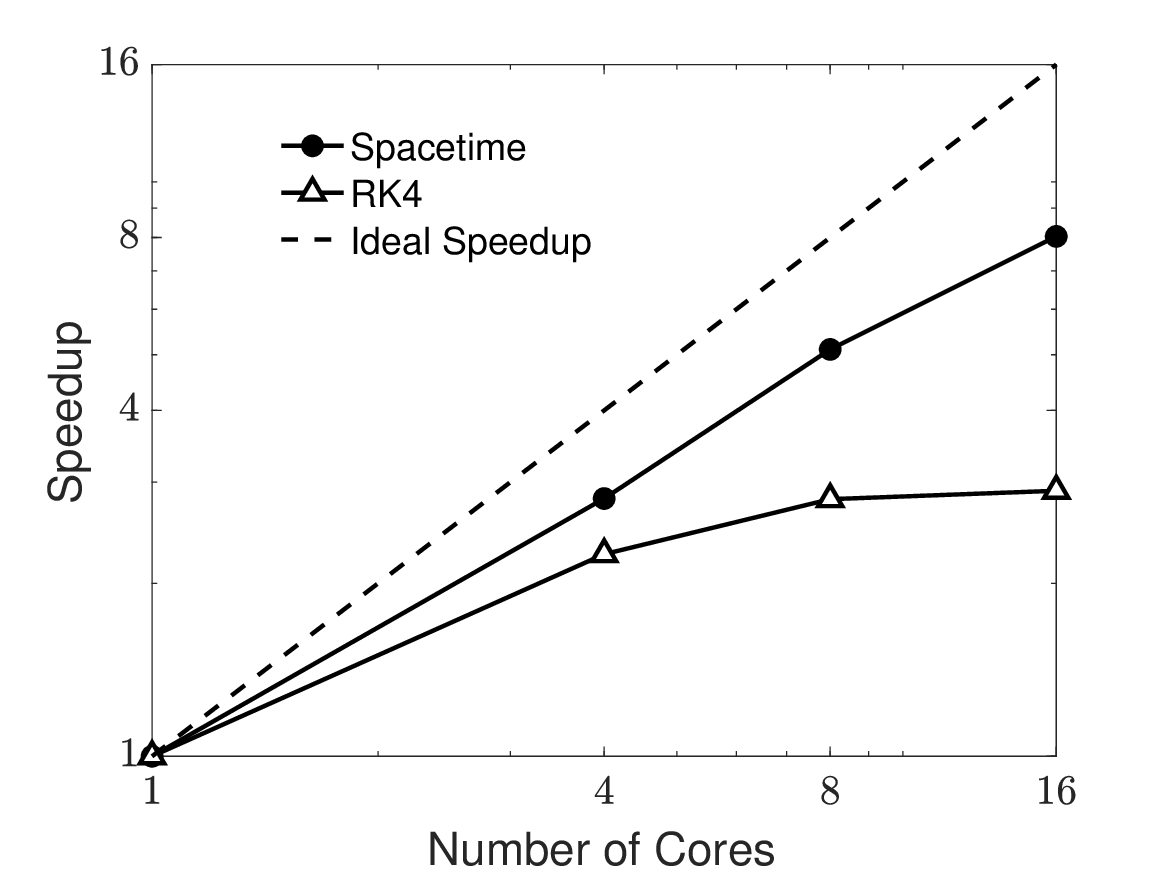}
    \caption{Strong scaling of the spacetime and RK4 methods for the shock advection problem at $j = 6$.}
    \label{fig:speedup6}
\end{figure}

\begin{figure}[H]
\begin{subfigure}{0.47\textwidth}
  \includegraphics[width=\linewidth]{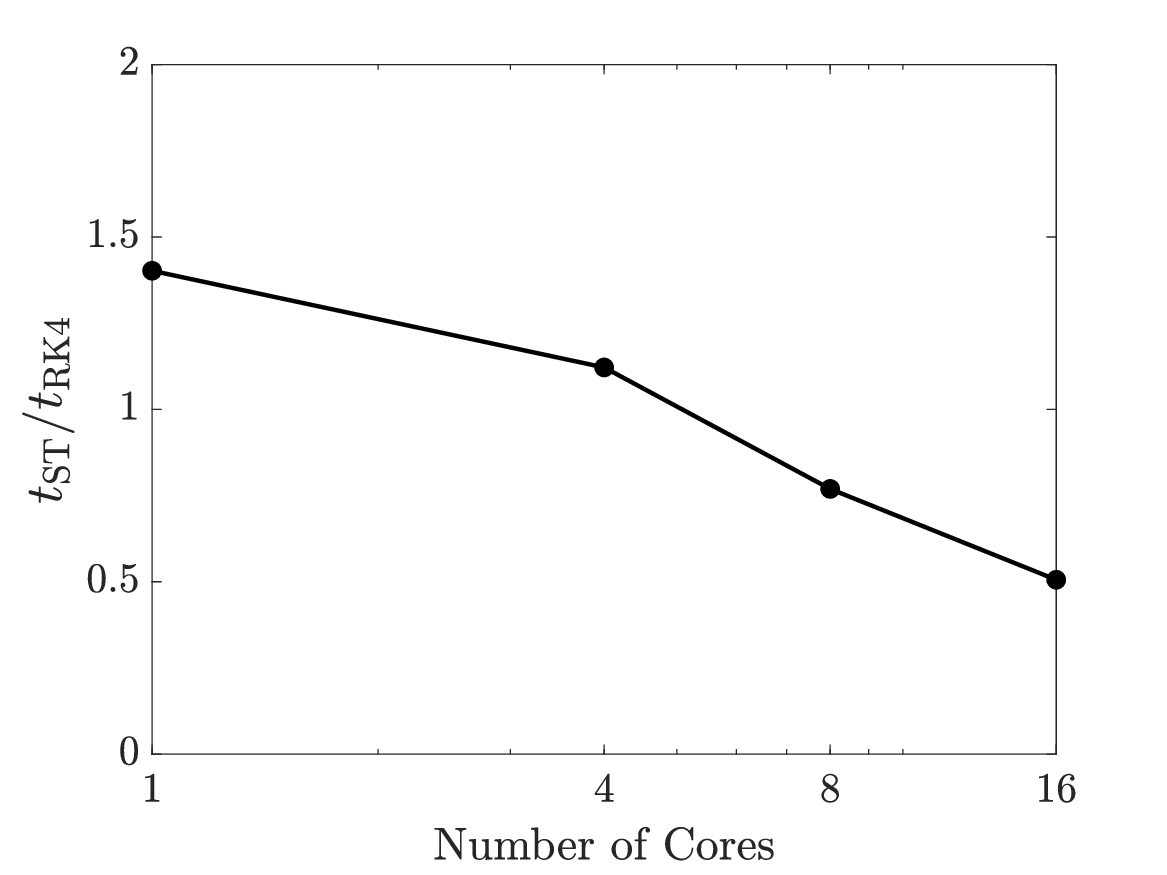}
  \caption{Relative time to solution of spacetime method vs. RK4 with varying number of cores at $j=6$.} \label{fig:rkDT}
\end{subfigure}
\hspace*{\fill}
\begin{subfigure}{0.49\textwidth}%
  \includegraphics[width=\linewidth]{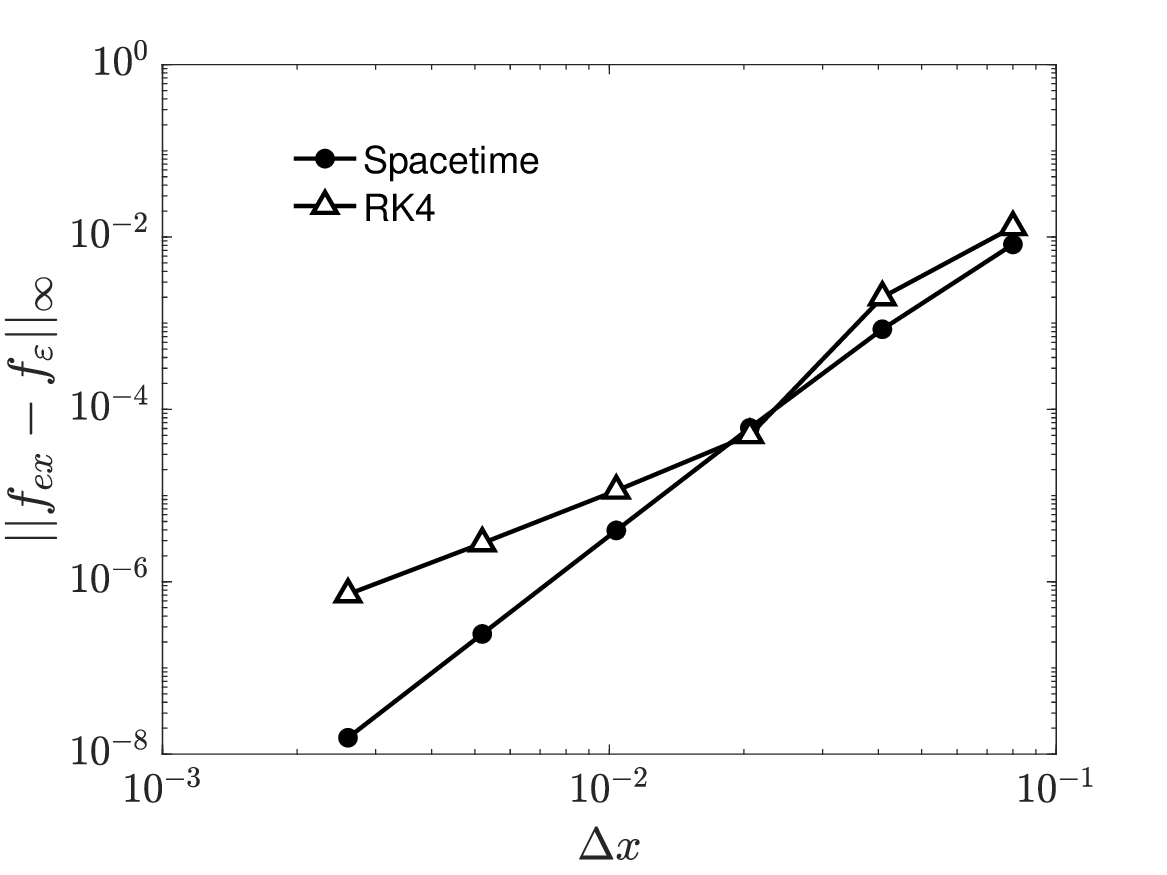}
  \caption{Solution convergence from $j = 1$ to $j = 6$ .}\label{fig:rkAcc}
\end{subfigure}
\caption{Spacetime wavelet vs. RK4 method time and accuracy comparisons for the shock advection problem.}
\label{fig:RKvsST}
\end{figure}

For these comparisons, we solve and analyze the shock advection problem (see Section \ref{sec:mb}) with the same parameters, spatial grid ($j = 6$), size of the domain, initial conditions, boundary conditions, and spatial basis order, $p_x = 6$. We use temporal basis order $p_t = 4$ for the spacetime results. We select CFL Courant number $C = 0.1$ based on traditional linear stability analysis which provides a stable solution and accuracy similar to the spacetime method with large $\Delta x$. The two methods are also constrained to output the same amount of data to ensure that a fair comparison is being made between them. The reported times in Figs. \ref{fig:speedup6} and \ref{fig:rkDT} are computed based on the average time to solution over five runs.

For this simple problem with relatively few degrees of freedom, the most time consuming block of the spacetime algorithm is building the tangent sparsity pattern and the computation of its entries. However, when properly parallelized with MPI, the construction time can be minimized. The MUMPS solver package also allows for the solve step of the Newton-Raphson loop to run in parallel. For these reasons, we achieve considerable speedup of the algorithm. The RK4 method is easily parallelized in space, but parallelizing in time, while possible \cite{gander2007analysis, gander201550, subber2016asynchronous, prakash2014computationally}, is a much more complicated process. Fig. \ref{fig:speedup6} shows the speedup of both algorithms when running on a varying number of cores.

We see that for the dense problem, the spacetime method is capable of achieving speedup closer to ideal and scales well with increasing number of cores. On the other hand, a time-marching method like RK4 does not benefit to the same degree when parallelized only in the spatial direction as typical for many solvers. For RK4 implementation in MRWT \cite{harnish2018adaptive}, the function updates required at each step are performed in parallel. However, each step and update in time is performed sequentially, as each step is dependent on the previous one \cite{verwer1996explicit, moin2010fundamentals}. The spacetime method achieves improved scaling due to the fact that the entire spacetime problem and the resulting system of algebraic equations are formulated in both space and time simultaneously, allowing for efficient parallel discretization and computation. This  leads to better parallel utilization as more collocation points (\textit{i.e.}, combining both spatial and temporal dimensions) are distributed across the computational node.  In contrast, the RK4 time-marching method stagnates early only at four cores, as it is parallelized only over spatial dimension. 
 
Next we make a comparison between the time to solution and accuracy for both methods. For the accuracy comparison, the viscosity parameter $\nu$ was relaxed from $0.01$ to $0.1$ to admit solutions for a wider range of $\Delta x$. We see from Fig. \ref{fig:rkDT} that the RK4 method outperforms the spacetime solver when run on fewer cores. However, the spacetime method surpasses RK4 when both are parallelized over higher number of cores. This again is due in part to the spacetime method allowing for natural parallelization in the time dimension. 
Fig. \ref{fig:rkAcc} illustrates the convergence rate of the solution for both the RK4 and the spacetime solvers. We note that the RK4 method offers $\mathcal{O}(\Delta t^4)$ truncation error, and at larger $\Delta x$ it maintains the overall $4$th-order convergence as the spatial error dominates. However,  at sufficiently small $\Delta x$ the CFL condition requires a small $\Delta t$. Suddenly, the overall $4$th-order convergence of RK4 solver is lost as the $\mathcal{O} (\sfrac{1}{\Delta t}$) accumulation error begins to dictate the convergence rate \cite{mosbach2009quantitative}. On the other hand, the spacetime method maintains  $4$th-order convergence for all $\Delta x$ considered. 
Note that for the RK4 method, an extremely small timestep can be chosen that may yield a more accurate solution and retain the convergence order. However, this small $\Delta t$ can handcuff the efficiency of the method by drastically increasing total time to solution \cite{de2013courant, verwer1996explicit}.

\subsubsection{Walking Burgers' Equation with Time-Dependent Advection (Shock Advection)} \label{sec:mbct}
\begin{figure}[H]
\begin{subfigure}{0.49\textwidth}
  \includegraphics[width=\linewidth]{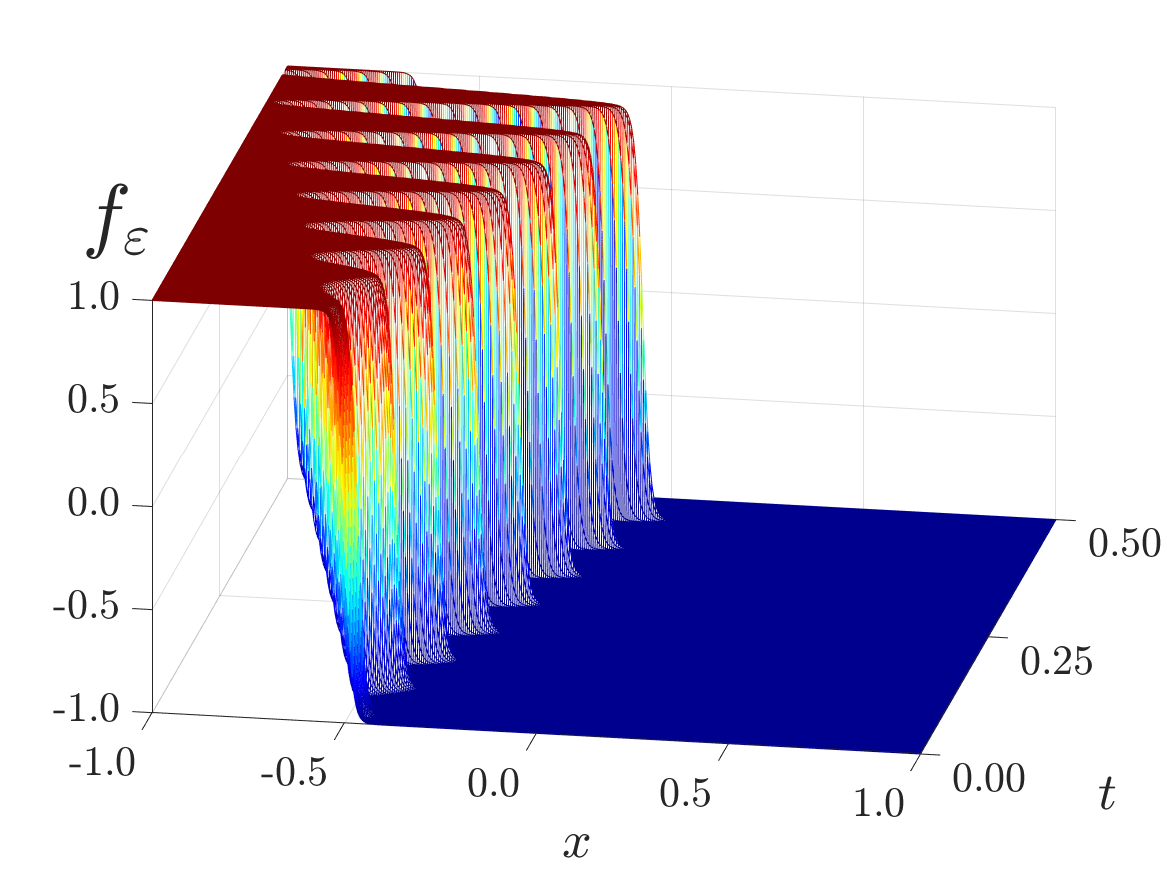}
  \caption{Spacetime solution.}\label{fig:mbctST}
\end{subfigure}
\hspace*{\fill}
\begin{subfigure}{0.49\textwidth}%
  \includegraphics[width=\linewidth]{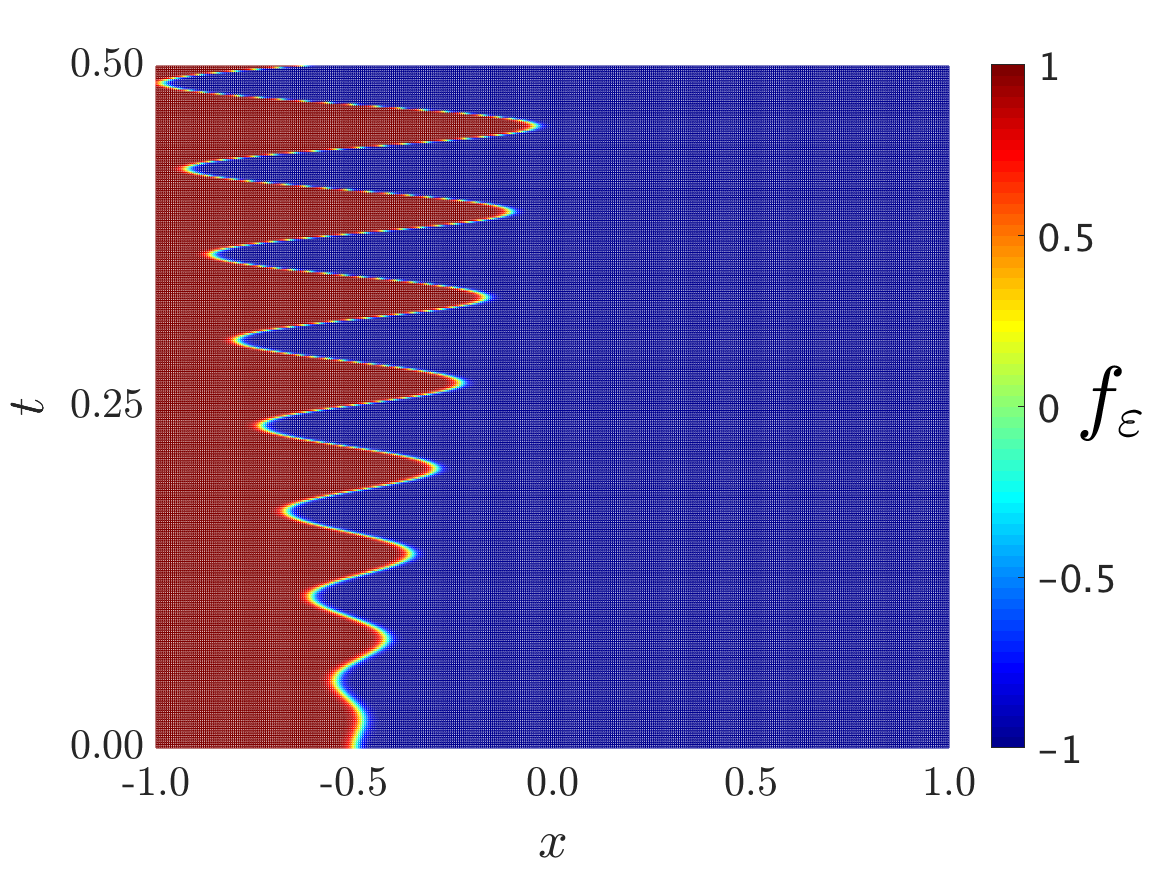}
  \caption{Top view.}\label{fig:mbctTop}
\end{subfigure}
\caption{Spacetime solution and profile for shock advection problem with $c(t) = \mathrm{sin}\left(\sfrac{t}{\tau}\right)$ at $j = 7$, $p_x =p_t= 6$.}
\label{fig:mbct3d}
\end{figure}

The capability of the solver to handle multiple scales in both the spatial and temporal directions can be observed by modifying the shock advection problem to use a time-dependent advection coefficient. The boundary and initial conditions remain identical to the $c(t) = 1$ case discussed above. Using the Method of Manufactured Solutions (MMS) \cite{salari2000code} to keep the system well-defined, we set $c(t) = \text{sin}\left(\sfrac{t}{\tau}\right)$, where $\tau = 10^{-2}$ is the chosen time scale of this modified problem. The spacetime solution of this example is shown in Fig. \ref{fig:mbct3d}, achieving accuracy of $\lVert f_{ex} - f_{\varepsilon}\rVert_{\infty} = 3.9 \times 10^{-6}$ at level $j = 7$ (2,362,369 DOF) with $p_x=p_t=6$. The thickness of the shock again is governed by the viscosity parameter $\nu$. However, the shock front advances in time sinusoidally due to the time-dependent forcing term introduced by MMS.

\begin{figure}[H]
\begin{subfigure}{0.49\textwidth}
  \includegraphics[width=\linewidth]{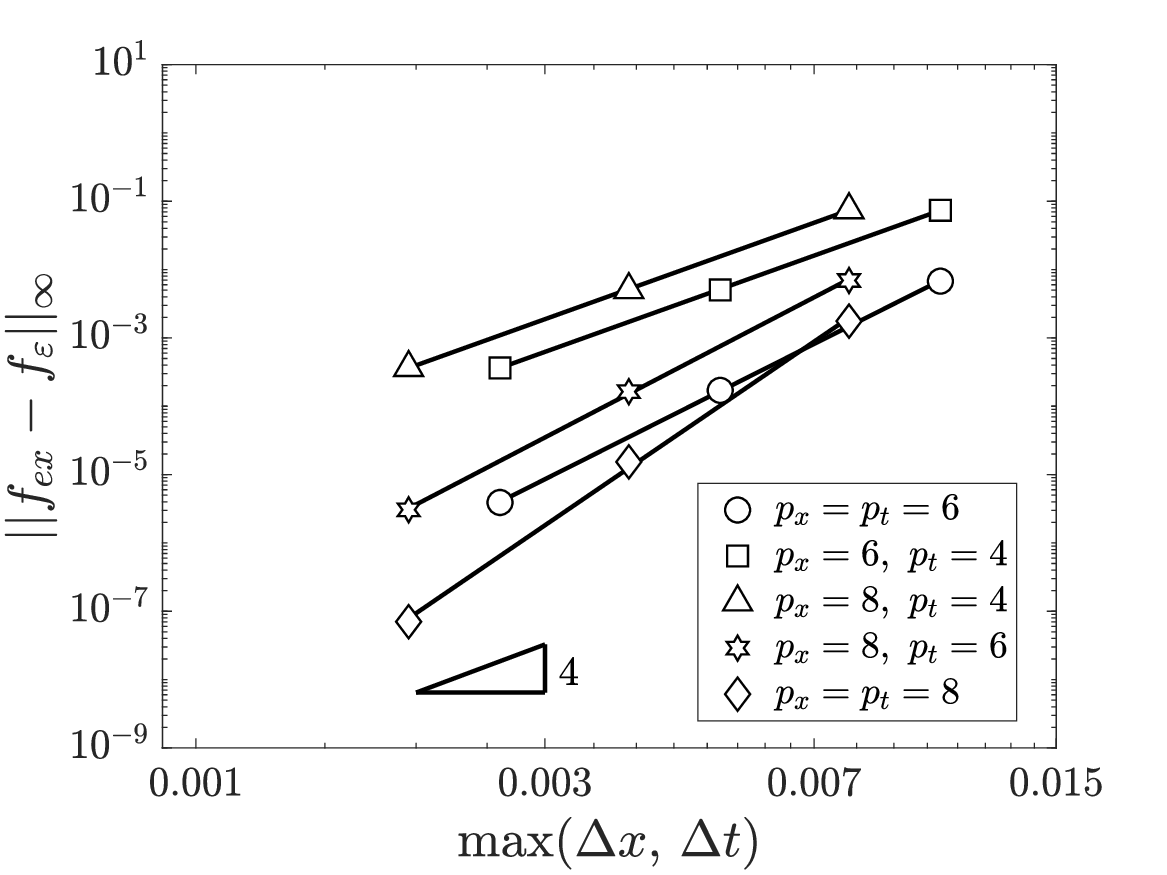}
  \caption{Solution convergence.}\label{fig:mbctSolConv}
\end{subfigure}
\hspace*{\fill}
\begin{subfigure}{0.49\textwidth}%
  \includegraphics[width=\linewidth]{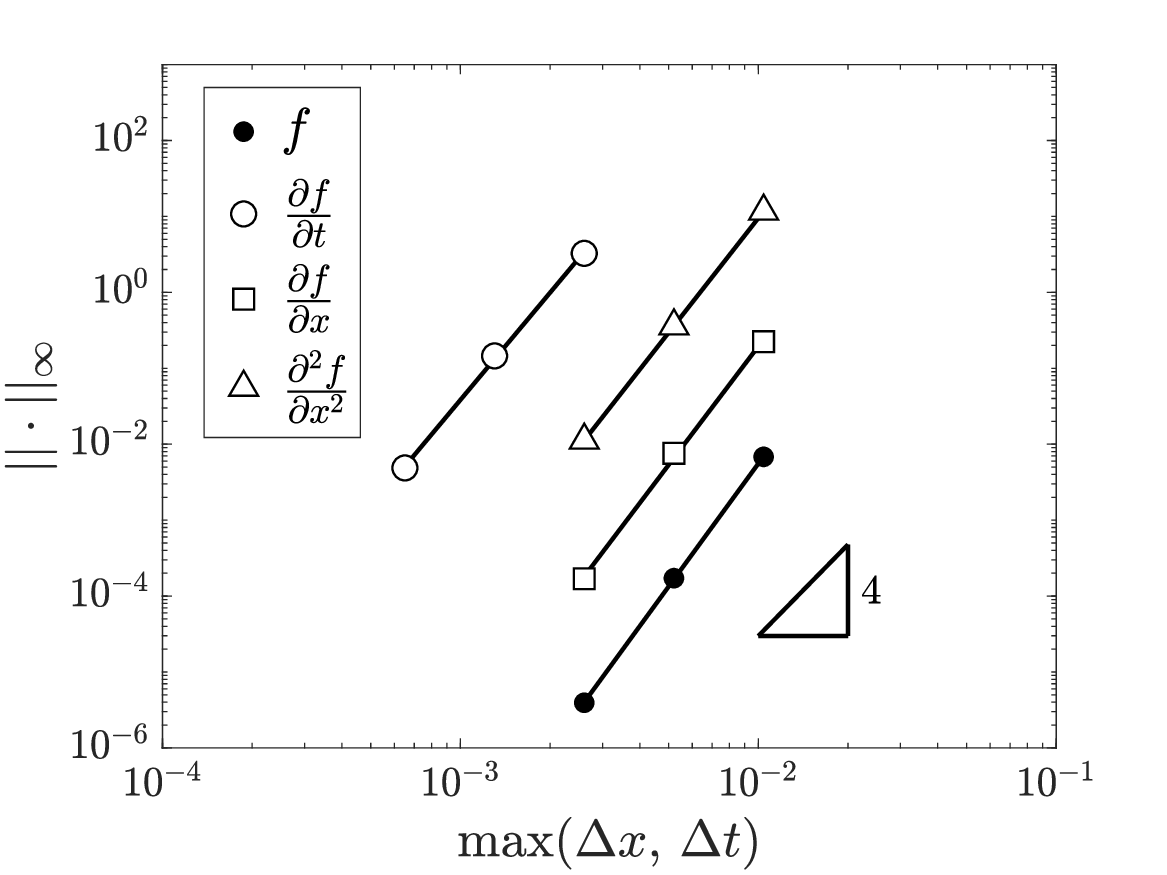}
  \caption{Solution and derivative convergence for $p_x =p_t=6$.}\label{fig:mbctDconv}
\end{subfigure}
\caption{Solution and derivative convergence for the shock advection problem with $c(t) = \mathrm{sin}\left(\sfrac{t}{\tau}\right)$ at $j = 5, 6, 7$.}
\label{fig:mbctConv}
\end{figure}

\noindent The convergence for both solution and derivative estimates are displayed in Fig. \ref{fig:mbctConv} in the same manner as in the $c(t) = 1$ example. 
Again, we achieve $4$th-order or better solution convergence for all combinations of basis functions (Fig. \ref{fig:mbctSolConv}), along with $4$th-order convergence for the solution and all derivative approximations as shown in Fig. \ref{fig:mbctDconv}. We note that this problem achieves better global accuracy by increasing the order of the temporal basis functions, due to its complex, time-varying solution. For example, $p_x = p_t = 6$ is more accurate than $p_x = 6$ and $p_t = 4$. This was not the case for the problem with $c(t) = 1$ (see Fig. \ref{fig:mbSolConv}), where both combinations achieved the same accuracy. The highest accuracy and convergence rate is obtained with $p_x=p_t=8$ ($\mathcal{O}(\Delta q^6)$). The capability of the spacetime wavelet solver to accurately solve this problem with high-order convergence is noteworthy due to the multiple regions with steep gradients in both space and time.

\subsubsection{Steepening Burgers' Equation (Shock Evolution)} \label{sec:steep}
\begin{figure}[H]
\begin{subfigure}{0.47\textwidth}
  \includegraphics[width=\linewidth]{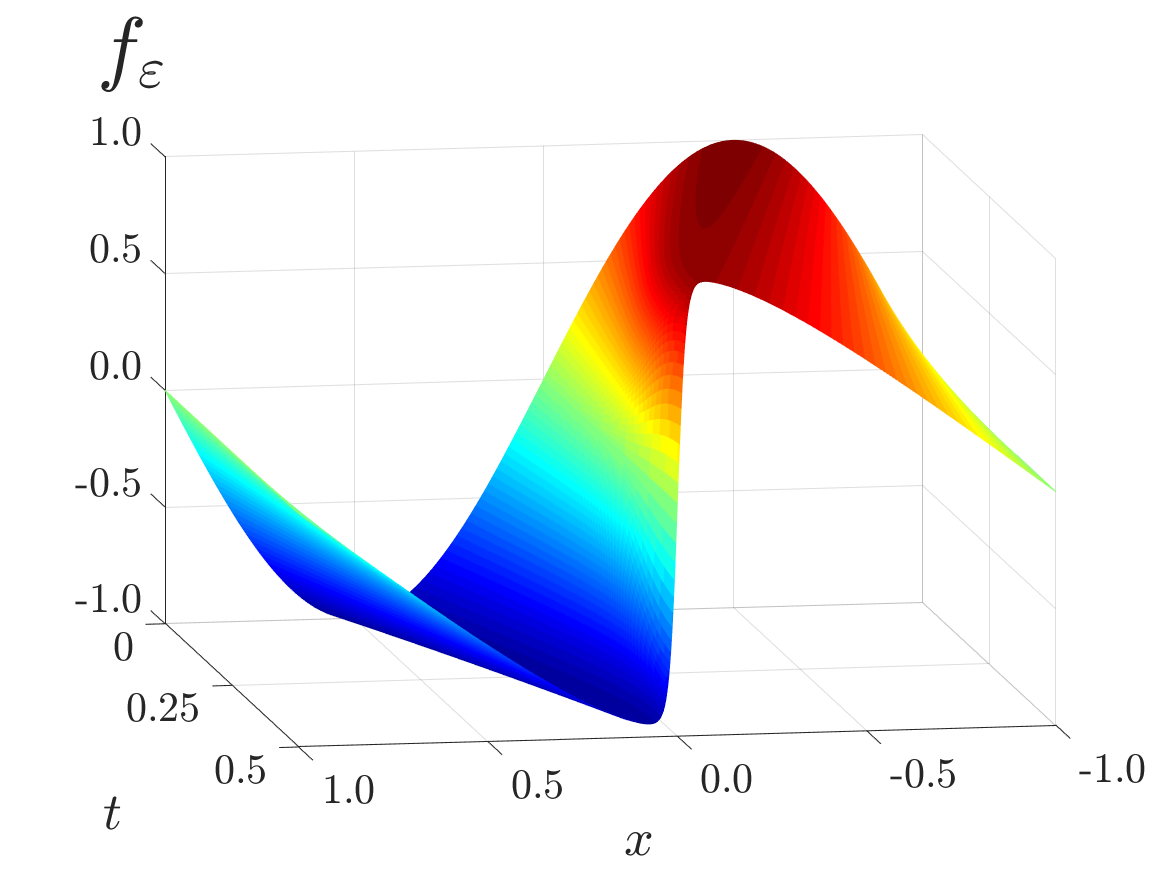}
  \caption{Shock evolution spacetime solution.}\label{fig:burgST}
\end{subfigure}
\hspace*{\fill}
\begin{subfigure}{0.49\textwidth}%
  \includegraphics[width=\linewidth]{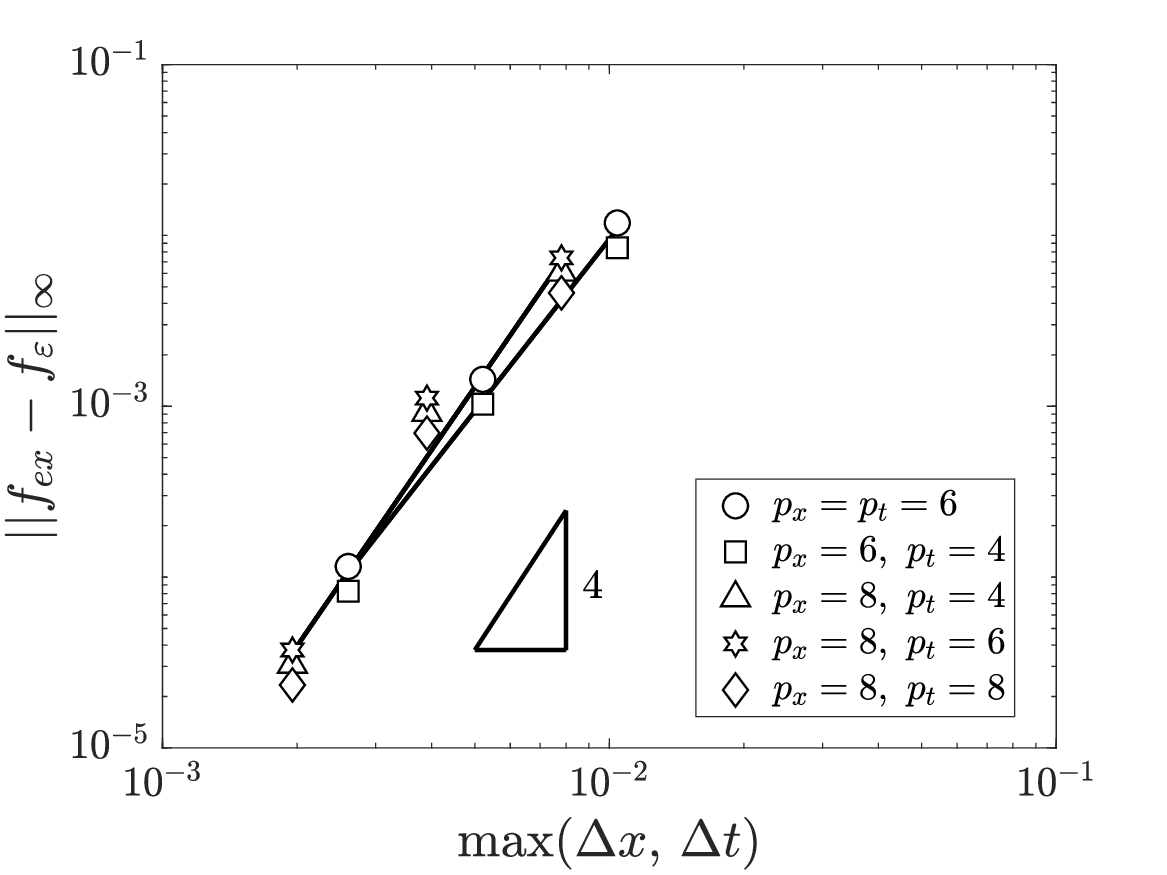}
  \caption{Shock evolution solution convergence.}\label{fig:burgSolConv}
\end{subfigure}
\caption{Spacetime solution with $p_x = 6$, $p_t = 4$ and convergence rates for the shock evolution problem for $j = 4, 5, 6$}
\label{fig:burg3d}
\end{figure}

The final Burgers' problem we solve is the shock evolution problem, defined by an advection coefficient $c(t) = 0$, eliminating the linear advection term (see Eq. (\ref{eq:modBurg})).
 Dirichlet boundary and initial conditions are defined by the analytical solution in the domain:
 \begin{align}
    f(x,t) &= - \frac{\int_{-\infty}^{\infty} \sin{( \pi x - \pi \eta )} \text{exp}\left( \frac{ - \cos{(\pi x - \pi \eta )} }{2 \pi \nu } \right) \text{exp}{\left( \frac{ - \eta^{2} }{ 4 \nu t } \right)} \mathrm{d} \eta }{\int_{-\infty}^{\infty} \text{exp}{\left( \frac{ - \cos{(\pi x - \pi \eta)}}{ 2 \pi \nu} \right)} \text{exp}{\left( \frac{ - \eta^{2} }{ 4 \nu t}\right)} \mathrm{d} \eta }.
    \label{eq:steepeningSol}
\end{align}
The spacetime solution and convergence rate of this problem are shown in Fig. \ref{fig:burg3d}, with relative accuracy $\lVert f_{ex} - f_{\varepsilon}\rVert_{\infty} = 6.7 \times 10^{-6}$  achieved at level $j = 6$ (394,497 DOF) with $p_x=6, p_t=4$.
\noindent We again see that $4$th-order or convergence or higher is achieved for all combinations of basis orders.

Despite using the same combinations of basis functions for all Burgers' solution convergence plots (Figs. \ref{fig:mbSolConv}, \ref{fig:mbctSolConv}, \ref{fig:burgSolConv}), we observe that the effect of the basis orders on solution accuracy is different. This is because the effect of the basis functions and their order in space and time along with solution features impact accuracy. For the shock problem with time-dependent advection (Fig. \ref{fig:mbctConv}), increasing $p_t$ increases the accuracy of the solution due to the rapidly-varying features in the temporal direction (see Fig. \ref{fig:mbct3d}). In contrast, for the shock steepening problem (Fig. \ref{fig:burg3d}), increasing the order of the basis functions in time has almost no effect on the solution error due to the relatively minimal changes in solution features throughout the time domain.

The results of these nonlinear Burgers' experiments demonstrate the capability of the spacetime wavelet method to accurately solve nonlinear systems with steep, shock-like features with high-order convergence rates as predicted by the theory.
\subsection{Sod Shock Tube} \label{sec:sst}
The Sod shock tube problem models a shock wave and a contact discontinuity moving to the right of the domain with a rarefaction wave moving back to the left \cite{sod1978survey}. The invicid Riemann solution to this problem features discontinuities. However this work uses hyperbolic tangent functions to mimic shock-like initial conditions while eliminating discontinuities by introducing viscosity. The Sod shock tube problem is governed by the 1D Navier-Stokes equations
\begin{align}
    \begin{split}
    &\frac{\partial \rho}{\partial t}  +   \frac{\partial(\rho v)}{\partial x}  =   0 \ \ \mathrm{in} \ \ \Omega,\\
    &\frac{\partial v}{\partial t}  +   \rho \nabla v \cdot v  -  \frac{\partial \sigma}{\partial x}  = \ 0 \ \ \mathrm{in} \ \ \Omega, \\
    &\rho \frac{\partial e}{\partial t}  +  \rho \nabla e \cdot v  +  \sigma : \nabla v  +  \text{div}\left( \hat{q} \right)  =  0 \ \ \mathrm{in} \ \ \Omega, \\
    &\rho = \rho_I, \ v= v_I, \ e = e_I \ \ \mathrm{on} \ \ \Omega_x \times (t = 0),  \\ 
     &\rho = \rho_B, \ v= v_B, \ e = e_B  \ \ \mathrm{on} \ \ \partial\Omega_x \times (0, T]. \label{eq:sod}
     \end{split}
\end{align}
Assuming a calorically perfect ideal gas, Fourier heat conduction and zero bulk viscosity, we obtain the constitutive equations
\begin{align}
    \hat{q} = \kappa \frac{\partial T}{\partial x}, \ \ \ \  e = c_v T, \ \ \ \ E = e + \frac{1}{2}v^2, \ \ \ \ \tau = \frac{4}{3}\mu \frac{\partial v}{\partial x},  \ \ \ \ \sigma = \tau - p,  \ \ \ \  p = (\gamma-1)\rho e. 
    \label{eq:sodConst}
\end{align}
The governing and constitutive equations are discretized with a dense wavelet representation similarly to the Burgers' examples, yielding the discrete residual equations
\begin{align}
\begin{split}
    \mathbf{R}_\rho  &=   \boldsymbol{\rho}\cdot{}^{(1, t)}\boldsymbol{\Gamma} +  \left({}^{(1, x)}\boldsymbol{\Gamma}\cdot\boldsymbol{\rho}\right) \circ \mathbf{v} + \boldsymbol{\rho} \circ \left({}^{(1, x)}\boldsymbol{\Gamma}\cdot\mathbf{v}\right),\\
    \mathbf{R}_v  &=  \boldsymbol{\rho}\circ \left[ \mathbf{v}\cdot\boldsymbol{\Gamma}^{(t)} + \left({}^{(1, x)}\boldsymbol{\Gamma}\cdot\mathbf{v}\right)\circ \mathbf{v}\right] + \left(\gamma-1\right)\left[\left({}^{(1, x)}\boldsymbol{\Gamma}\cdot\mathbf{v}\right)\circ \mathbf{e} + \boldsymbol{\rho} \circ \left({}^{(1, x)}\boldsymbol{\Gamma}\cdot\mathbf{e}\right) \right] - \frac{4}{3}\mu \left({}^{(2, x)}\boldsymbol{\Gamma}\cdot\mathbf{v}\right),\\
    \mathbf{R}_e  &=  \boldsymbol{\rho}\circ \left[ \mathbf{e}\cdot{}^{(1, t)}\boldsymbol{\Gamma} + \mathbf{v}\circ \left({}^{(1, x)}\boldsymbol{\Gamma}\cdot\mathbf{e}\right)\right] + \left[\frac{4}{3}\mu \left({}^{(1, x)}\boldsymbol{\Gamma}\cdot\mathbf{v}\right) - \left(\gamma-1\right)\boldsymbol{\rho} \circ \mathbf{e} \right]\circ \left({}^{(1, x)}\boldsymbol{\Gamma}\cdot\mathbf{v}\right) - \frac{\kappa}{c_v}\left({}^{(2,x)}\boldsymbol{\Gamma}\cdot\mathbf{e}\right).
    \end{split}
\end{align}
For the sake of brevity, we will not list the nine equations describing the tangent matrix $\mathbf{K}$, but they are defined in the Appendix \ref{sec:sodTan}. The system in the  $\mathbf{A}\vec{x} = \vec{b}$ form is visualized as
\begin{align}
\begin{bmatrix}
\mathbf{K}_{\rho\rho} & \mathbf{K}_{\rho v} & \mathbf{0}\\
\mathbf{K}_{v\rho} & \mathbf{K}_{vv} & \mathbf{K}_{ve} \\
\mathbf{K}_{e\rho} & \mathbf{K}_{ev} & \mathbf{K}_{ee}
\end{bmatrix}
\begin{bmatrix}
\vec{\rho} \\ 
\vec{v} \\ 
\vec{e}
\end{bmatrix}
\ = \ -
\begin{bmatrix}
\vec{\mathcal{R}}_\rho \\ 
\vec{\mathcal{R}}_v \\ 
\vec{\mathcal{R}}_e
\end{bmatrix},
\label{eq:sodkMat}
\end{align}
where $\mathbf{K}_{mn}$ is the tangent block resulting from the derivative of $\mathcal{R}_m$ with respect to unknown  $n$.
The material parameters used are consistent with those of dry air at room temperature and are given in Table \ref{tab:parameters}.
\begin{table}[H]
\begin{center}
\renewcommand{\arraystretch}{2}
\begin{tabular}{ |c|c|c| } 
\hline
Variable & Name & Value \\ 
 \hline\hline
 $\gamma$ & Ratio of specific heats & $7/5$ \\ 
 \hline
 $\kappa$ & Thermal conductivity & $2.55 \times 10^{-2} \ \mathrm{W}/(\mathrm{m} \cdot \mathrm{K})$ \\ 
 \hline
 $c_v$ & Constant volume specific heat &  $7.18 \times 10^2 \ \mathrm{J}/(\mathrm{kg} \cdot \mathrm{K})$ \\ 
 \hline
 $\mu$ & Dynamic viscosity & $1.9 \times 10^{-5}  \ \mathrm{Pa} \cdot \mathrm{s}$ \\
 \hline
\end{tabular}
\end{center}
\caption{Material parameters of dry air at room temperature.}
\label{tab:parameters}
\end{table}
We use the boundary conditions defined by Kamm \cite{kamm2008enhanced} on the domain $x \in [0, 1] \ \text{cm}$, $t \in [0, 0.2] \ \text{s}$, displayed in Table \ref{tab:bc}. The shock transition region is centered about $x = 0.5 \ \text{cm}$. Due to the stiffness in the initial conditions, the solution of the Sod problem cannot be obtained directly at the highest resolution level $j_{\mathrm{max}}$. Theretofore,  we utilize our novel recursive wavelet algorithm presented in Section \ref{sec:recursive} (see Algorithm \ref{alg:recAlg}).
\begin{table}[H]
\begin{center}
\renewcommand{\arraystretch}{2}
\begin{tabular}{ |c|c|c| } 
\hline
Density $\rho$ [$\text{g}/\text{cm}^3$] & Velocity $v$ [$\text{cm}/\text{s}$]& Pressure $p$ [$\text{dyn}/\text{cm}^2$]\\ 
 \hline\hline
 $\rho (0, t) = 1.0$  & $v(0,t) = 0.0$ & $p(0,t) = 1.0$ \\ 
 \hline
 $\rho (1, t) = 0.125$  & $v(1,t) = 0.0$ & $p(1,t) = 0.1$ \\
 \hline
\end{tabular}
\end{center}
\caption{Sod shock Dirichlet boundary conditions.}
\label{tab:bc}
\end{table}
\noindent To illustrate the significance of the weighting vectors in the recursive solution process we look at the initial condition for density, expressed as a hyperbolic tangent of the form
\begin{equation}
    \rho(x, t = 0) = A + B \ \mathrm{tanh}\left(\frac{x-x_0}{\delta}\right),
\end{equation}
where $A$ and $B$ are constant values and $\delta$ determines the steepness of the shock. The smaller the $\delta$ value, the steeper the curve. We gradually steepen the initial condition with the weighting vector $\vec{\delta} = [5 \delta, \ 3\delta, \ \delta]$. For this problem, the boundary condition forcing term has minimal impact on the system of equations, so the boundary condition weighting vector $\vec{\chi}$ contains values that enforce the complete set of boundary conditions for each solve as discussed in Section \ref{sec:NR}.

\begin{figure}[H]
    \centering
    \includegraphics[width=0.7\linewidth]{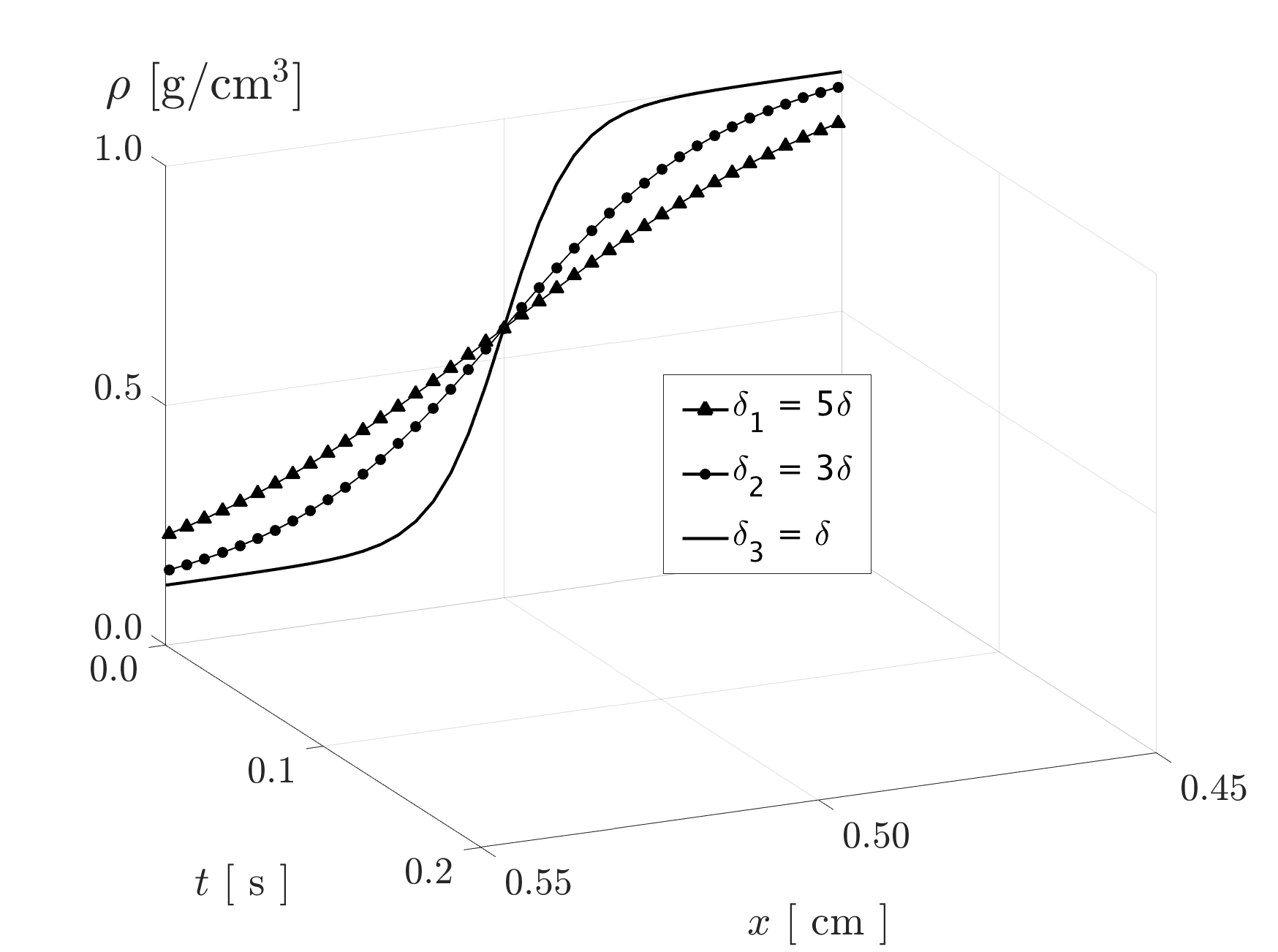}
    \caption{Recursive progression of Sod initial conditions resulting from Algorithm \ref{alg:recAlg}.}
    \label{fig:ICs}
\end{figure}

\begin{figure}[H]
\begin{subfigure}{0.49\textwidth}
  \includegraphics[width=\linewidth]{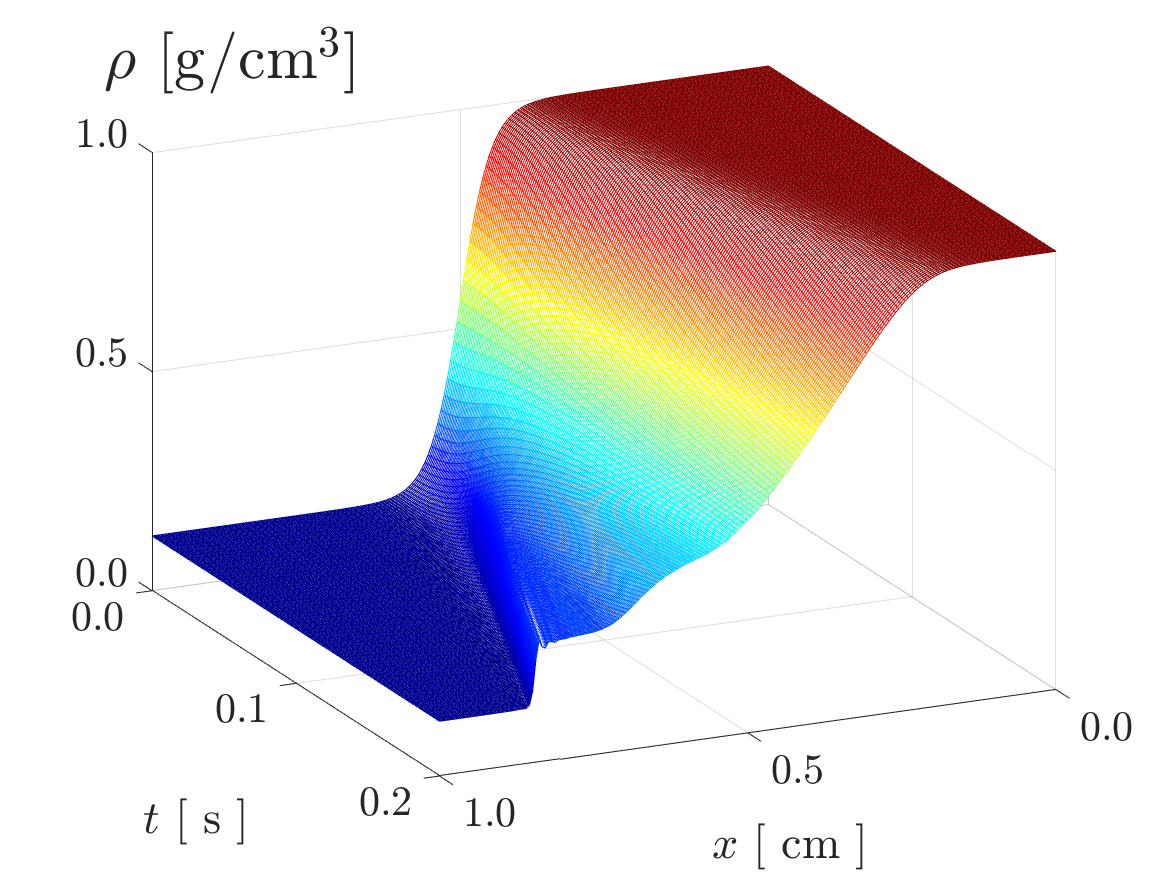}
  \caption{Sod density solution using $j=4, \delta_1 = 5 \delta$.}\label{fig:rhoj4}
\end{subfigure}
\hspace*{\fill}
\begin{subfigure}{0.49\textwidth}%
  \includegraphics[width=\linewidth]{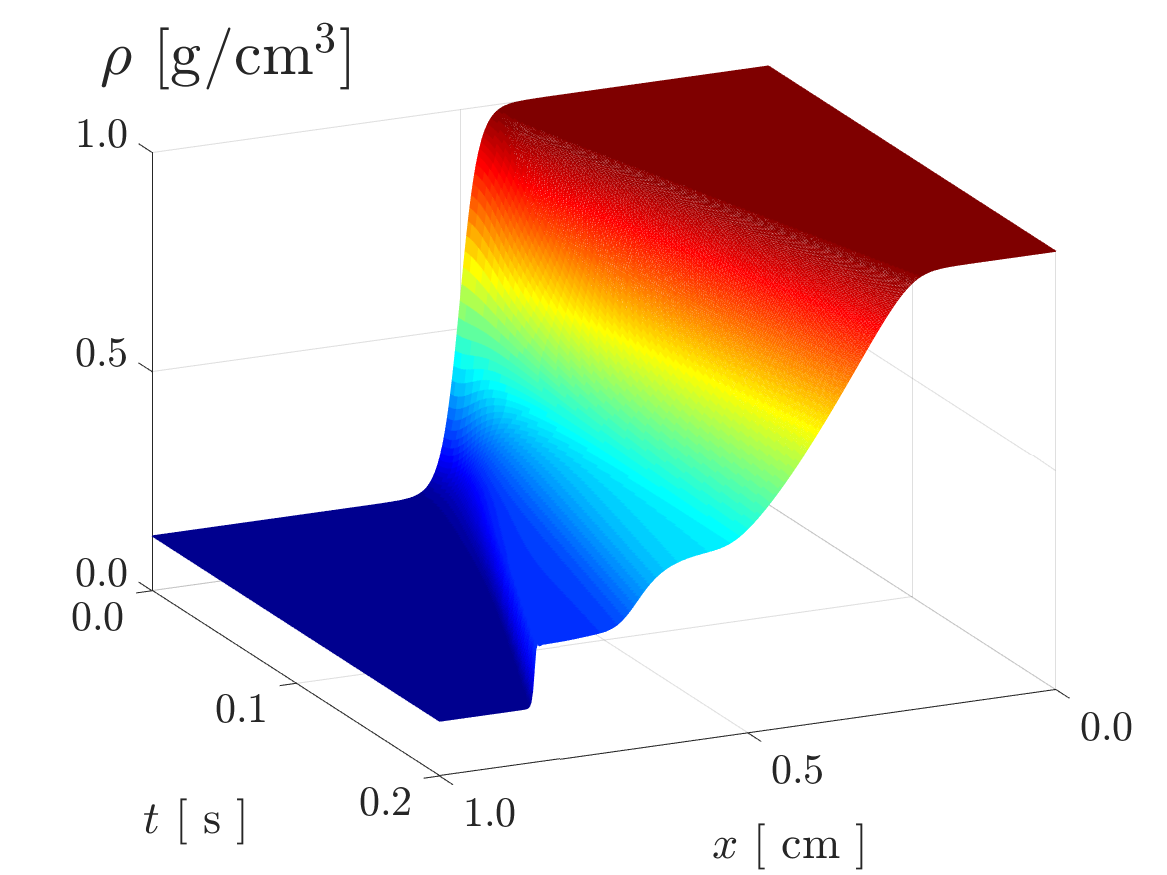}
  \caption{Sod density solution using $j=5, \delta_2 = 3 \delta$.}\label{fig:rhoj5}
\end{subfigure}
\caption{Sod solutions using recursive technique (Algorithm \ref{alg:recAlg}) at $j = 4$ and $j = 5$.}
\label{fig:rhoRec}
\end{figure}

The progression of steepening initial conditions used to solve for density with $\delta = 0.01$, $\mathrm{A} = \sfrac{9}{16}$, and $\mathrm{B} = -\sfrac{7}{16}$ is shown in Fig. \ref{fig:ICs}.
Moving from $j=4$  with $\delta_1 = 5 \delta$ to $j=5$ with $\delta_2
= 3 \delta$, we recursively reduce the sensitivity of the system,
Eq. (\ref{eq:solvability}), at each level by using a wavelet-generated
initial guess and relaxed initial conditions. These progressive
spacetime solutions are shown in Fig. \ref{fig:rhoRec}.

\noindent With an accurate solution obtained at $j=5$, we are now able to synthesize an initial guess that reduces system sensitivity at $j=6$ and satisfies the solvability condition, Eq. (\ref{eq:solvability}), allowing us to solve the Sod problem using the desired steep initial condition (\textit{i.e.}, $\delta = 0.01$). 

\begin{figure}[H]
\begin{subfigure}{0.49\textwidth}
  \includegraphics[width=\linewidth]{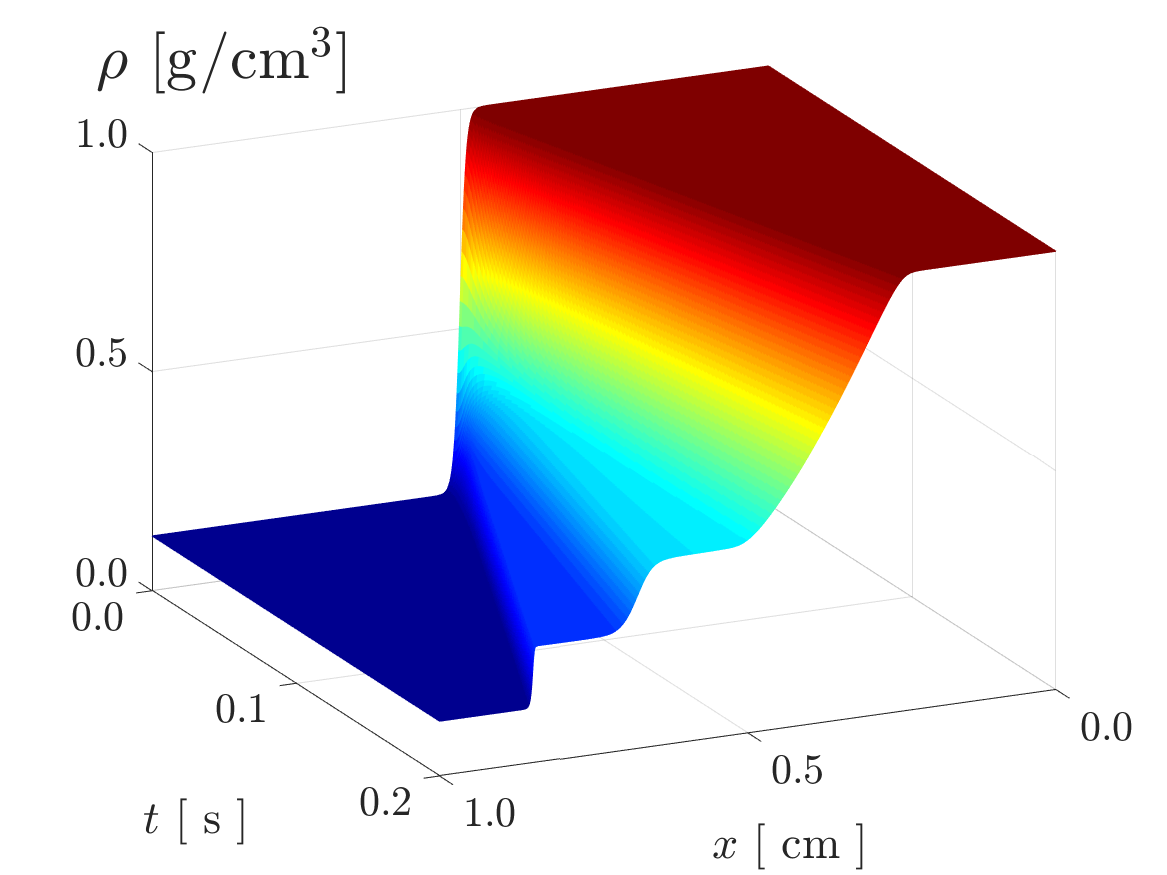}
  \caption{Density solution using spacetime wavelet solver.}\label{fig:rhoST}
\end{subfigure}
\hspace*{\fill}
\begin{subfigure}{0.49\textwidth}%
  \includegraphics[width=\linewidth]{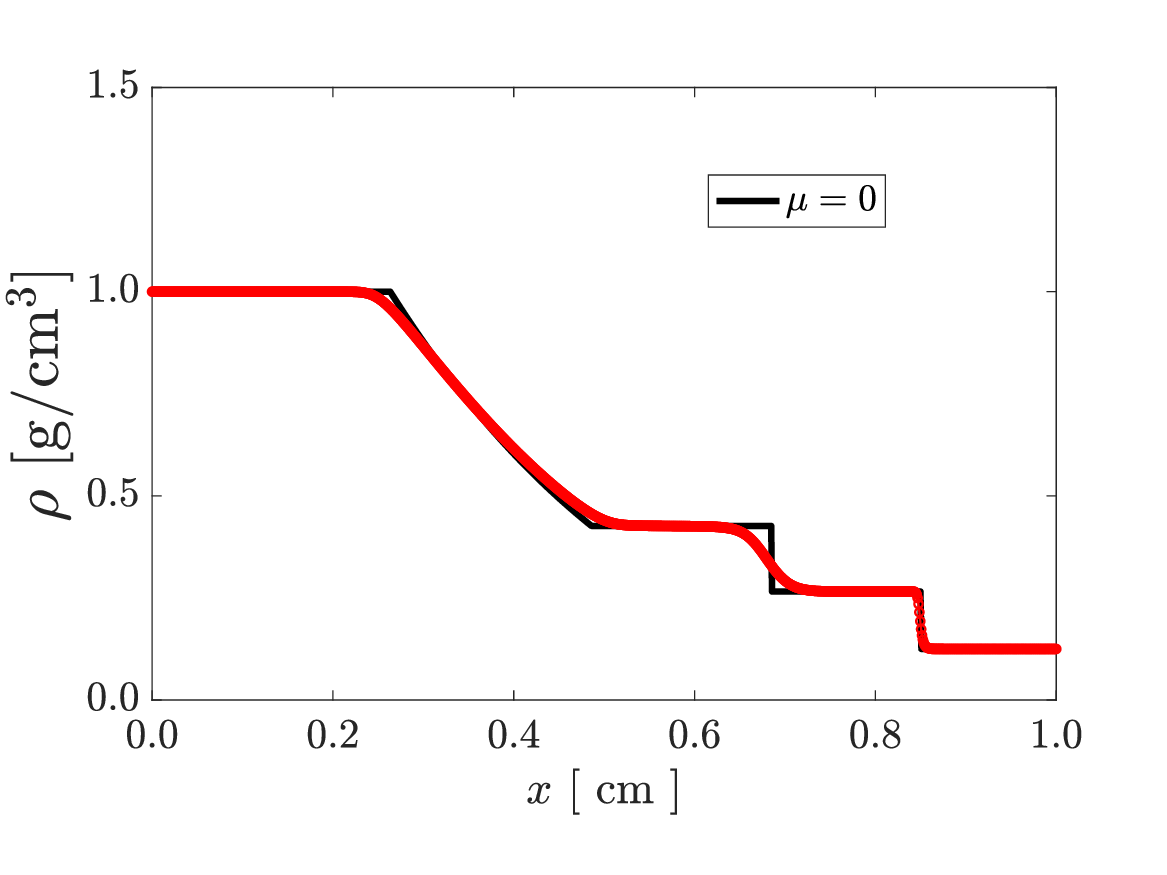}
  \caption{Density solution at $t = 0.2$ s plotted along with invicid solution.}\label{fig:rho2}
\end{subfigure}
\caption{Sod density solution using the spacetime solver at $j = 6$.}
\label{fig:rho}
\end{figure}

\begin{figure}[H]
\begin{subfigure}{0.49\textwidth}
  \includegraphics[width=\linewidth]{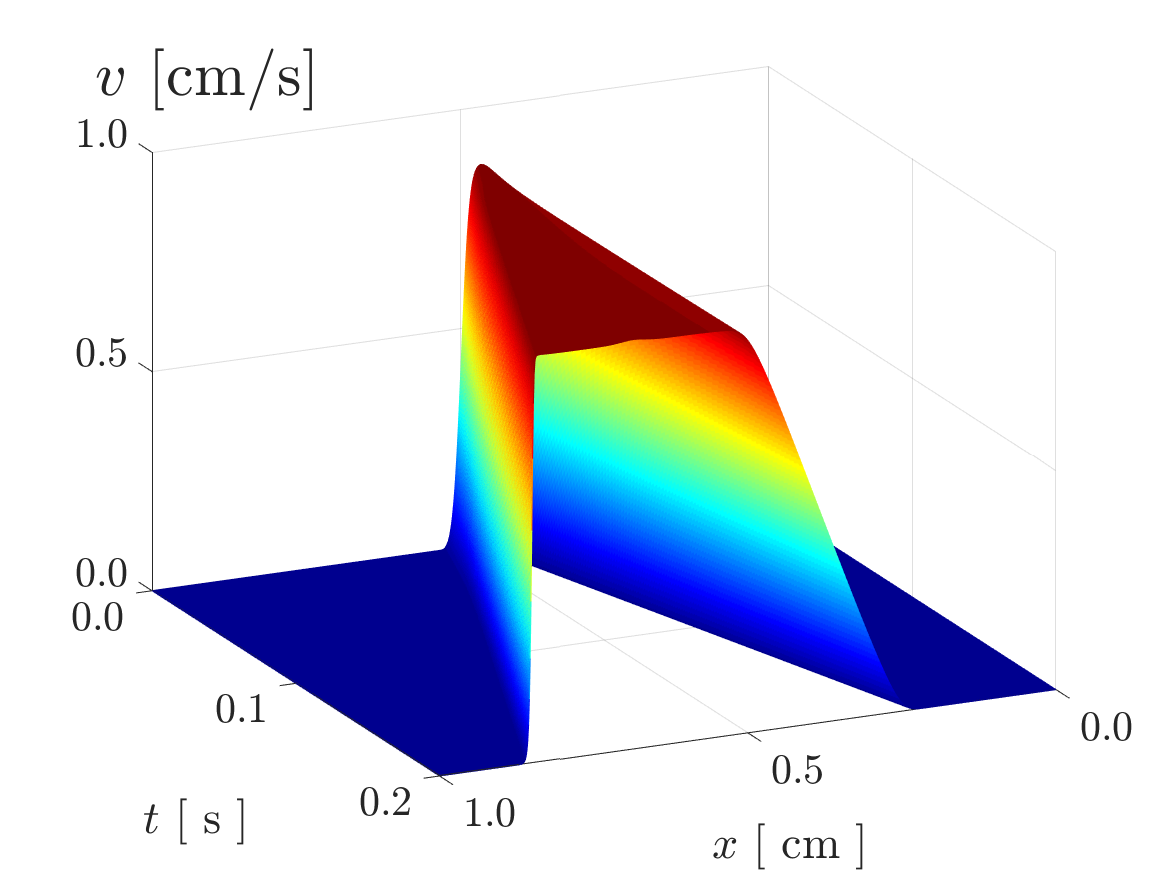}
  \caption{Velocity solution using spacetime wavelet solver.}\label{fig:vST}
\end{subfigure}
\hspace*{\fill}
\begin{subfigure}{0.49\textwidth}%
  \includegraphics[width=\linewidth]{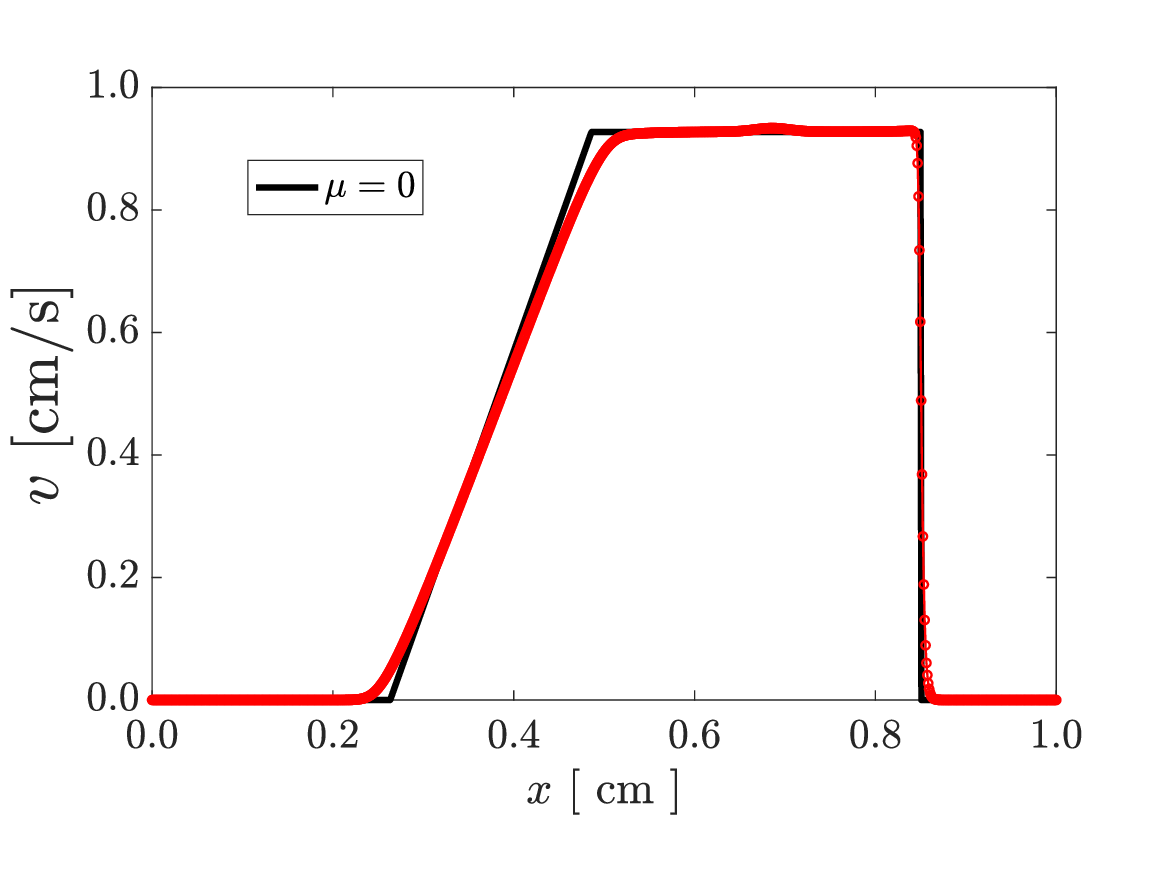}
  \caption{Velocity solution at $t = 0.2$ s plotted along with invicid solution.}\label{fig:v2}
\end{subfigure}
\caption{Sod velocity solution using the spacetime solver at $j = 6$.}
\label{fig:v}
\end{figure}

\begin{figure}[H]
\begin{subfigure}{0.49\textwidth}
  \includegraphics[width=\linewidth]{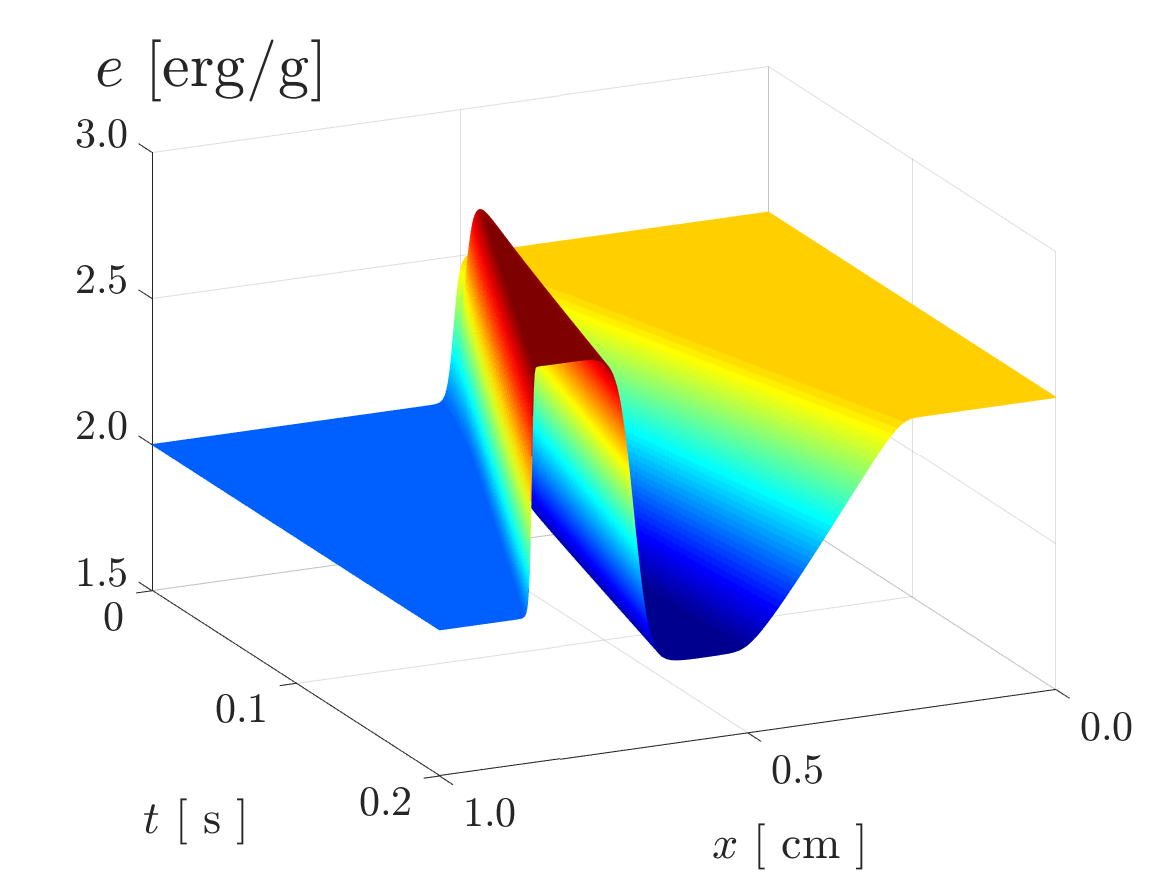}
  \caption{Energy solution using spacetime wavelet solver.}\label{fig:eST}
\end{subfigure}
\hspace*{\fill}
\begin{subfigure}{0.49\textwidth}%
  \includegraphics[width=\linewidth]{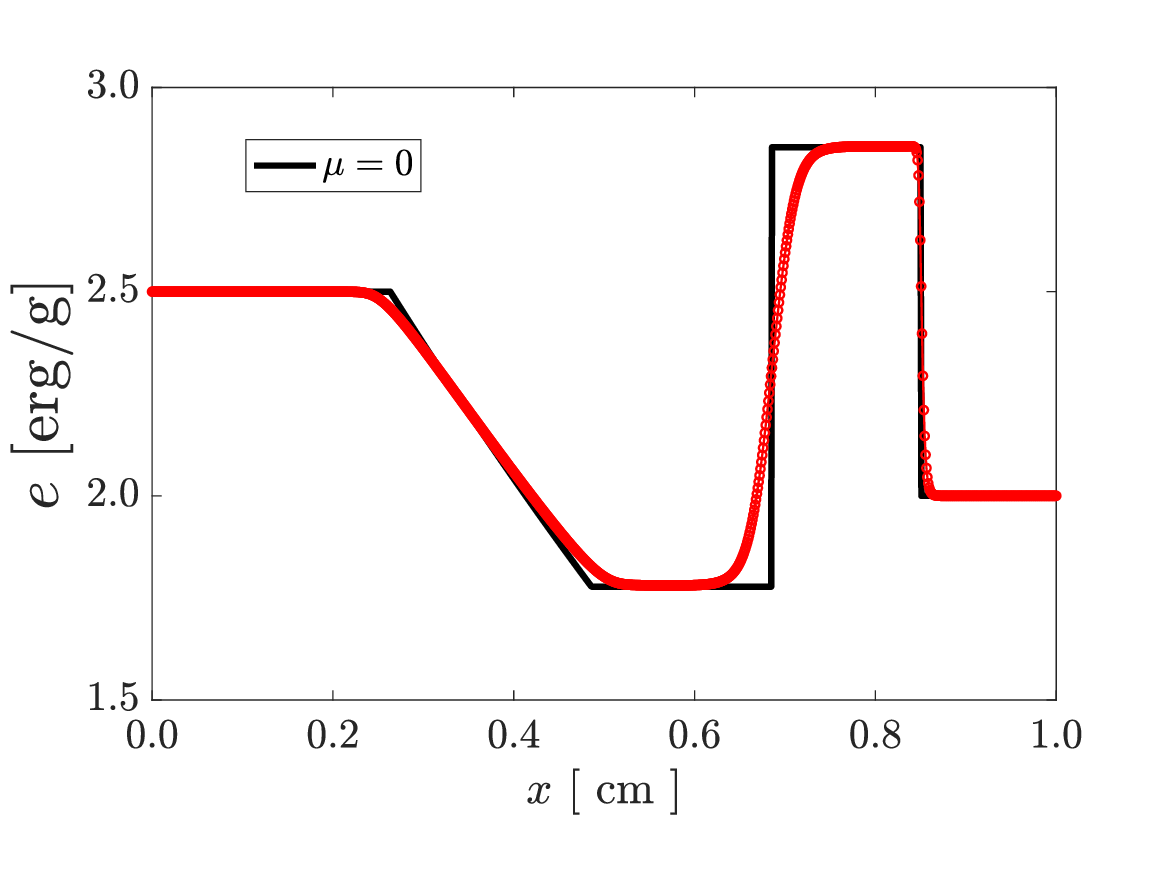}
  \caption{Energy solution at $t = 0.2$ s plotted along with invicid solution.}\label{fig:e2}
\end{subfigure}
\caption{Sod energy solution using the spacetime solver at $j = 6$.}
\label{fig:e}
\end{figure}

The final spacetime solutions with 1,774,083 total DOF are displayed in Figs. \ref{fig:rho}-\ref{fig:e}. Moreover, we present comparisons to the Riemann inviscid solutions at the final time $t = 0.2$ s. 
We see that the spacetime wavelet solver captures the rapidly-changing features in both the spatial and temporal directions. The steep viscous shock can be seen around  $x = 0.85 \ \mathrm{cm}$ in Fig. \ref{fig:v}. The contact discontinuity and refraction wave are well developed and can be observed in density profile (see Fig. \ref{fig:rho}). Transforming the solutions to wavelet space and analyzing wavelet coefficients at level $j_{\mathrm{max}}+1$, we use the relationship 
\begin{align}
    \mathrm{error}(\bullet) \propto \lVert \mathbb{d}_{\vec{k}}^{j_{\mathrm{max}+1}}\rVert_{\infty}
\end{align}
 to estimate $0.17\%$ error for density, $0.57\%$ error for velocity, and $1.36\%$ error for energy \cite{harnish2023adaptive}. These solutions were obtained using basis orders $p_x = p_t = 6$. Using a $6$th-order basis in time provides more temporal grid points and is helpful when capturing the rapidly-evolving features of the Sod solution (see also problem \ref{sec:mbct}). These results are obtained using a globally-convergent Newton-Raphson method, which keeps the initial solve with the uninformed guess from diverging.
 
\begin{figure}[H]
\begin{subfigure}{0.49\textwidth}
  \includegraphics[width=\linewidth]{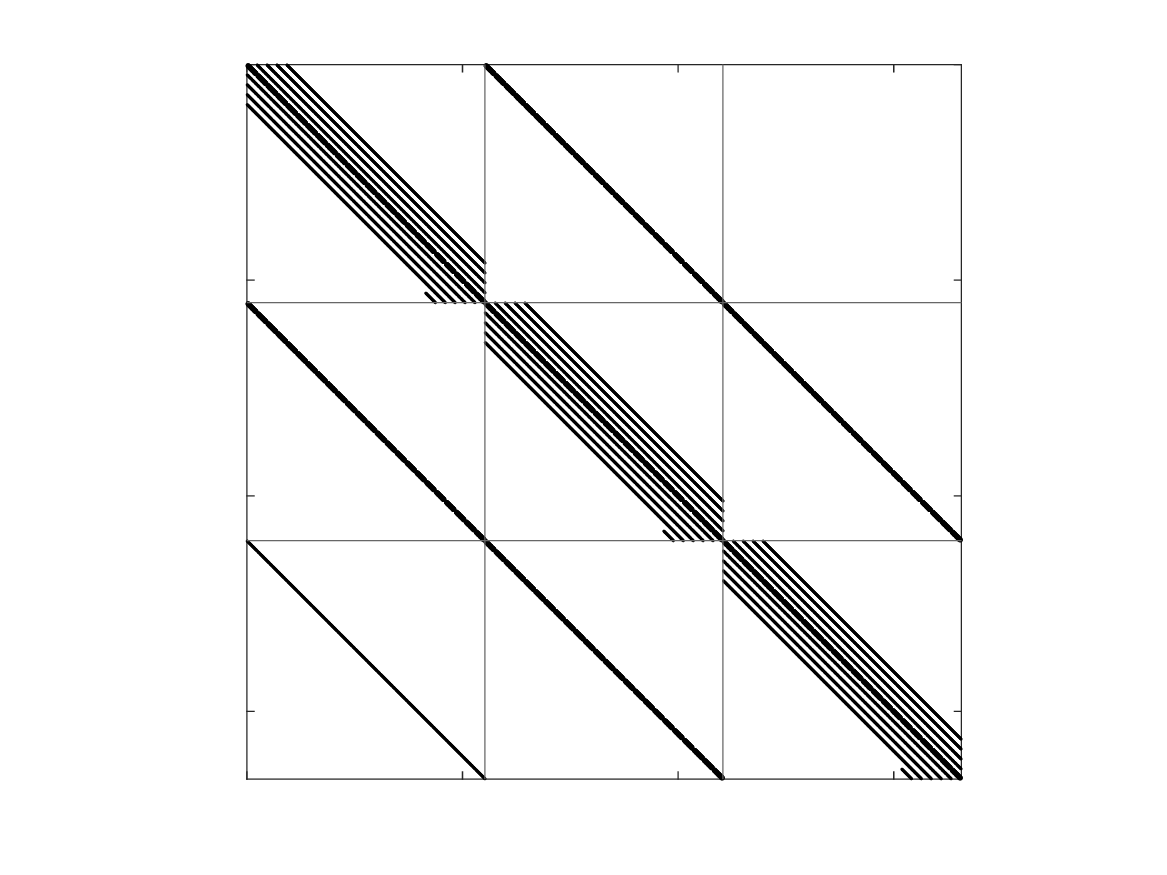}
  \caption{Sod problem sparsity pattern}\label{fig:sodSparsity}
\end{subfigure}
\hspace*{\fill}
\begin{subfigure}{0.49\textwidth}%
  \includegraphics[width=\linewidth]{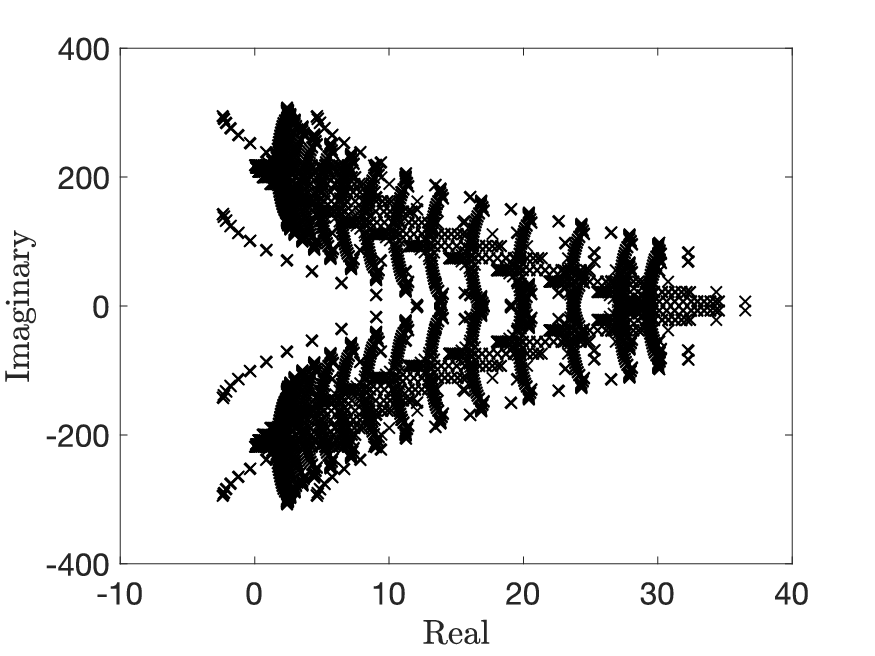}
  \caption{Sod problem eigenvalue spectrum}\label{fig:sodEigs}
\end{subfigure}
\caption{Spacetime Sod problem tangent sparsity pattern and eigenvalue spectrum at $j = 3$ with $p_x=p_t=6$.}
\label{fig:spEig}
\end{figure}

The effect of multiple governing equations on the structure and conditioning of the system can be observed in the sparsity pattern and eigenvalue spectrum of the tangent matrix described in Eq. (\ref{eq:sodkMat}). These patterns are displayed in Fig. \ref{fig:spEig}, and similarly to the discretized Burgers' problems (see Fig. \ref{fig:mbEigs}), we see that the spacetime Sod problem has entirely complex eigenvalues. Note that these plots are made by analyzing the tangent at level $j = 3$ with shallower initial conditions ($\delta = 0.1$) in order to more clearly show relevant features.

In Fig. \ref{fig:sodSparsity}, we see the nine distinct blocks, representing each of the Jacobian submatrices given in Eq. (\ref{eq:sodkMat}) (see Appendix \ref{sec:sodTan} for full list of equations). We see that the $K_{\rho\rho}$, $K_{vv}$ and $K_{ee}$ diagonal blocks of the matrix have similar structures to the Burgers' Jacobian (Fig. \ref{fig:mbSparsity}) for the same reasons listed previously. Again, the bandwidth of main diagonal bands and the number of off-diagonal bands are determined by problem features and the orders of the basis functions. 

\subsubsection{Effect of Recursive Algorithm on Convergence}
\begin{figure}[H]
    \centering
    \includegraphics[width=0.6\textwidth]{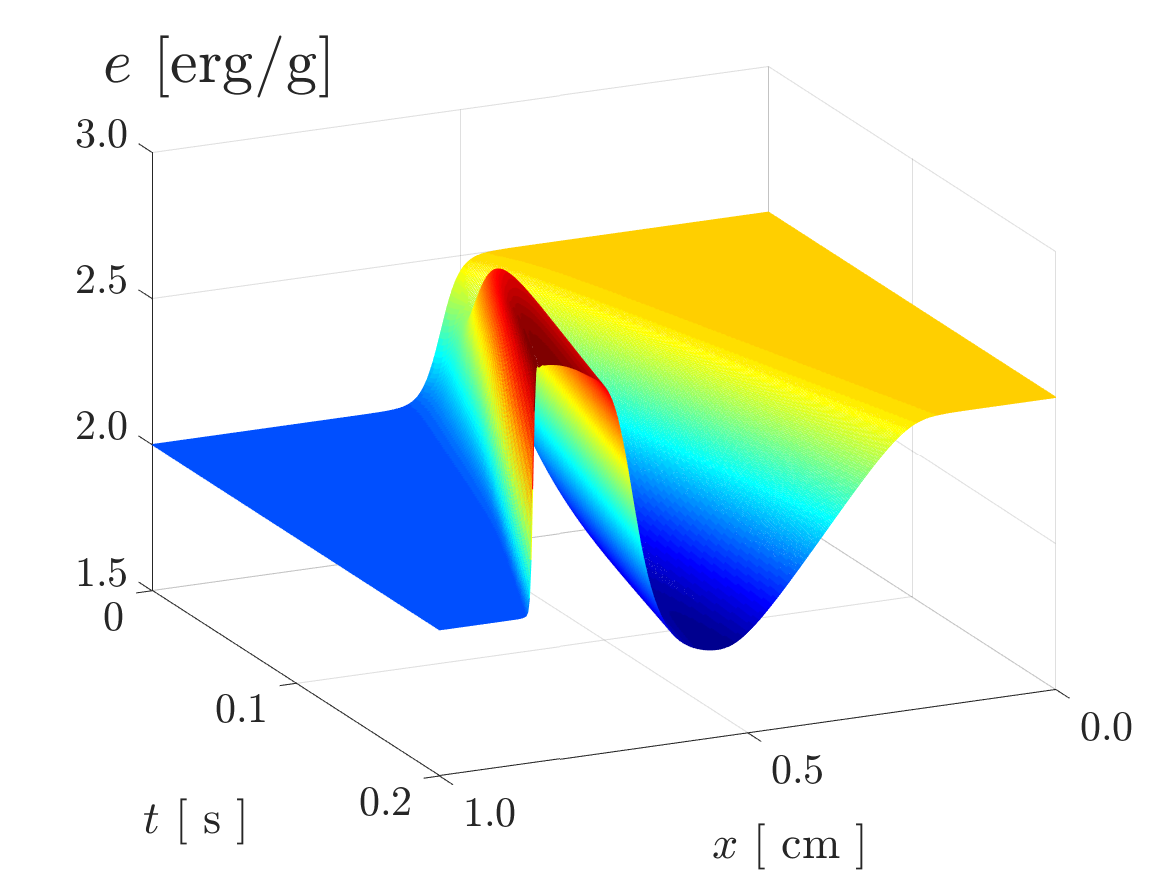}
    \caption{Spacetime energy solution at $j=5$ with relaxed viscosity $\delta = 0.05$.}
    \label{fig:ej5}
\end{figure}

\begin{figure}[H]
\begin{subfigure}{0.49\textwidth}
  \includegraphics[width=\linewidth]{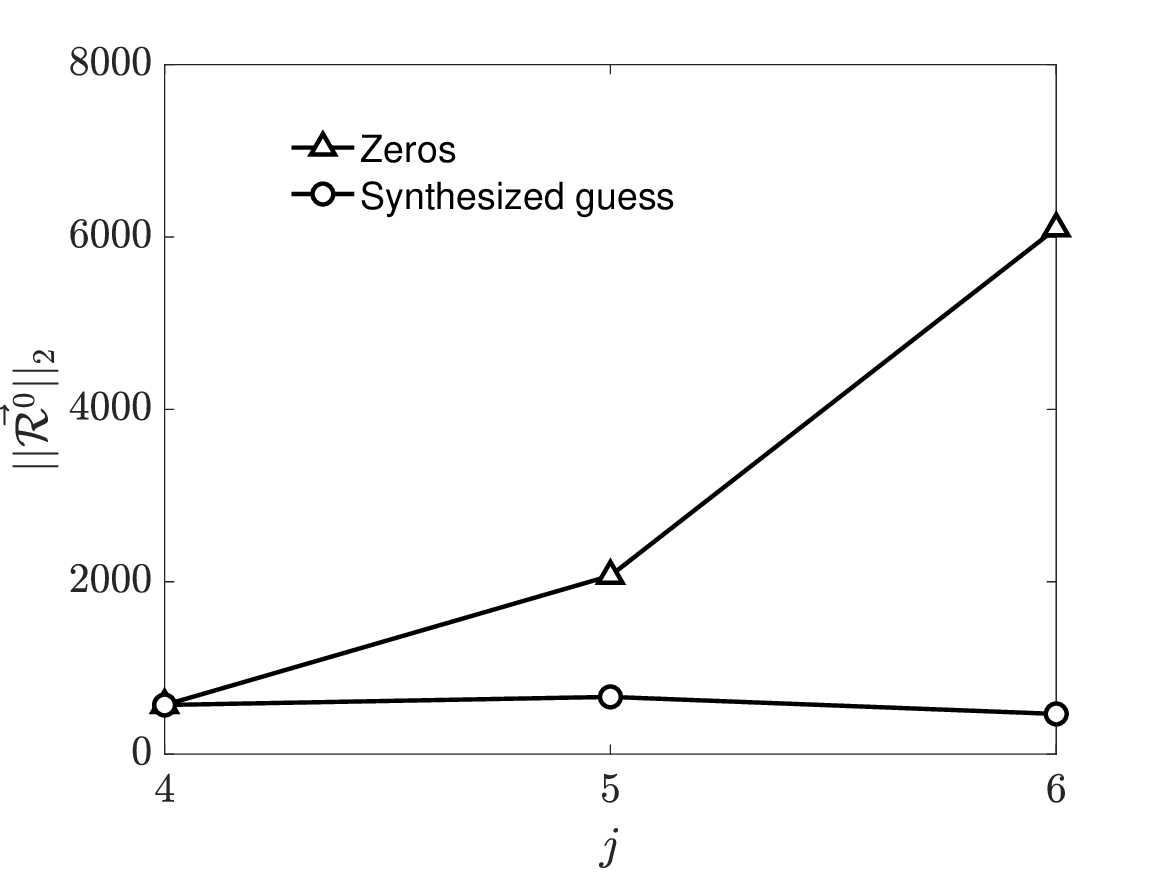}
  \caption{Initial right-hand side magnitude vs. resolution level}\label{fig:rVj}
\end{subfigure}
\hspace*{\fill}
\begin{subfigure}{0.49\textwidth}%
  \includegraphics[width=\linewidth]{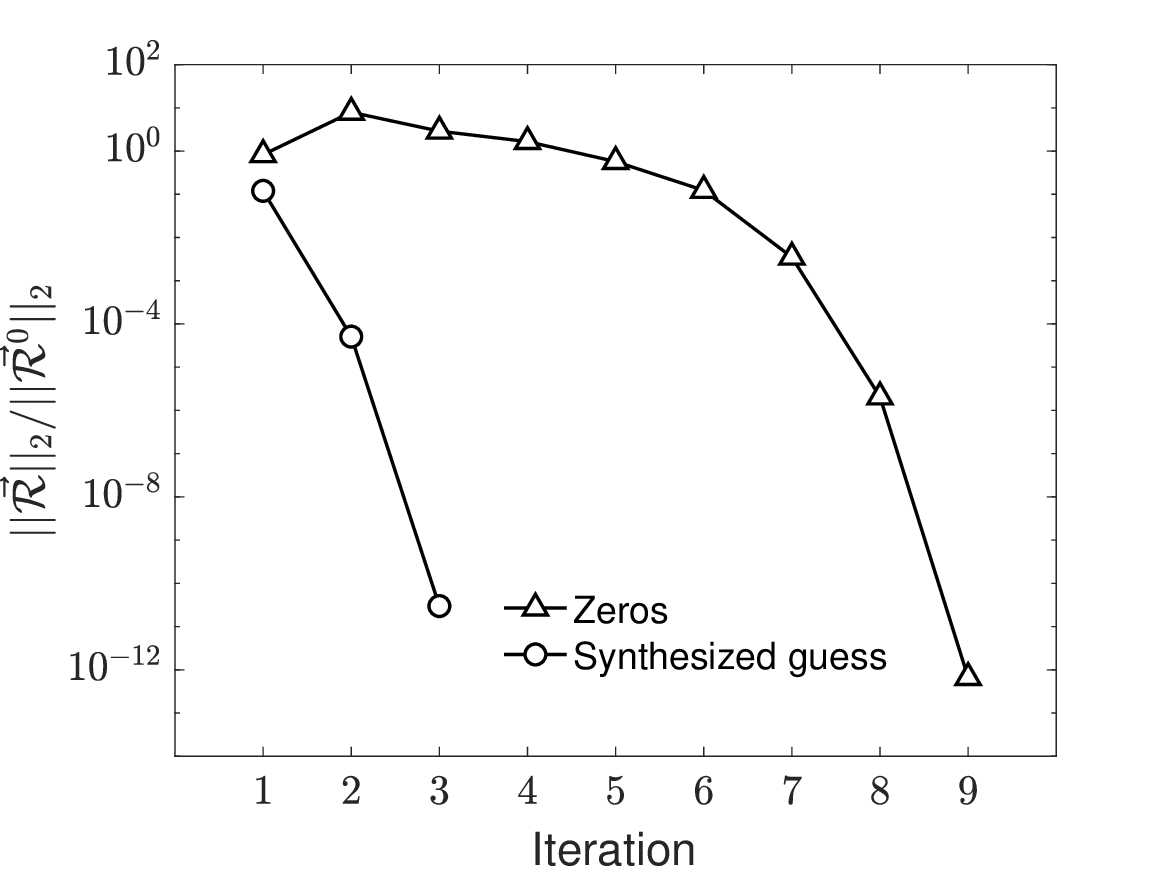}
  \caption{Relative error of Newton-Raphson method at $j=5$.}\label{fig:j5conv}
\end{subfigure}
\caption{Effect of wavelet synthesized initial guess on initial residual magnitude and solution convergence.}
\label{fig:recConv}
\end{figure}

For the full Sod problem at $j=6$ with steepness parameter $\delta =
0.01$, the recursive algorithm is necessary to converge to the desired
solution. Using an uninformed initial guess will cause the
Newton-Raphson method to begin outside the region of convergence and
ultimately diverge. In order to analyze the effect of the recursive
algorithm on the convergence of a system, we solve a relaxed problem
with $\delta = 0.05$. The relaxed spacetime solution for energy at
$j=5$ (148,225 DOF) is shown in Fig. \ref{fig:ej5}. As one can see,
the features are less sharp compared to full solution depicted in
Fig. \ref{fig:eST}. 
 Fig. \ref{fig:rVj} shows the magnitude of the initial residual vector
 using both a standard method (\textit{i.e.}, we solve directly at
 $j_{\mathrm{max}}$ without recursion) with a guess of zeros and the
 wavelet-based recursive technique using wavelet synthesis across
 $j=4, 5,$ and $6$. Fig. \ref{fig:j5conv} displays a Newton-Raphson
 convergence comparison at $j=5$ with a recursive initial guess from
 level 4 and with a zero initial guess.

We see that despite obtaining convergence with both guesses with this
relaxed problem, the synthesized guess reduces the magnitude of the
initial residual vector, Fig \ref{fig:rVj}, and moves the problem
closer to the region of convergence, resulting in faster convergence
(\textit{i.e.}, in 3 Newton steps compared to 9). It is worth noting
that even as the problem size grows substantially as resolution level
increases, the magnitude of the initial residual vector from the
wavelet-generated initial guess remains nearly
constant. Fig. \ref{fig:recConv} shows that when using an initial
guess of zeros, the residual initially grows but eventually decreases,
converging towards the solution. This is not always the case, and in
many problems, the residual will continue to grow and the solution
will diverge due to an overly-sensitive system. This problem is
avoided by using Algorithm \ref{alg:recAlg} to generate initial
guesses and reduce the sensitivity of the system,
Eq. (\ref{eq:solvability}). We note that for the Sod problem, the
wavelet initial guesses do not have a significant impact on the
conditioning of the tangent matrix.
\section{Conclusions} \label{sec:conclusions}
In this work, we have developed a novel spacetime high-order wavelet
method for problems with multiple spatial and temporal scales that
develop during shock loading, for example. We have demonstrated the
capability of a spacetime wavelet solver to obtain accurate solutions
to nonlinear PDE systems modeling shock behavior with high-order
global convergence and accuracy. This novel algorithm utilizes wavelet
approximation to discretize both the spatial and temporal components
of PDEs, including derivative terms, with \textit{a priori} error
estimates. We have developed a wavelet-based recursive algorithm to
improve system sensitivity and generate an optimal initial guess for
the Newton-Raphson method that is used to linearize the system. The
spacetime solver alleviates some of the concerns encountered when
utilizing time-marching methods, including timestep selection and
stability conditions. The numerical implementation is performed using
the Multiresolution Wavelet Toolkit, and the system of linear
algebraic equations is solved using the MUltifrontal Massively
Parallel sparse direct Solver. The wavelet theory describes high-order
convergence rates for both solution and derivative estimates, which
have been rigorously verified using several types of Burgers'
equation. Finally, we have demonstrated the spacetime method to handle
the coupled system of nonlinear Navier-Stokes equations with the Sod
shock tube results. Future studies will focus on the development of
the sparse spacetime representation. 
\section*{Acknowledgments}
This work was supported by the Network, Cyber, and Computational
Sciences Branch of the Army Research Office (ARO) under Award Number
W911NF-24-1-0039. Dr. Radhakrishnan Balu served as program monitor. We
would like to also thank Dr. Cale Harnish and Dr. Luke Dalessandro for
their assistance through discussions on wavelet theory and
computational implementation.
\newpage

\bibliography{newBib}

\begin{appendices}
\section{Sod Jacobian Equations}\label{sec:sodTan}
\begin{align*}
\begin{split}
\frac{\partial \mathbf{R}_\rho}{\partial \boldsymbol{\rho}} = \mathbf{K}_{\rho\rho} =& \ \mathbf{I}\otimes{}^{(1, t)}\boldsymbol{\Gamma} \ + \ \Bigl({}^{(1, x)}\boldsymbol{\Gamma} \otimes \mathbf{I}\Bigl) \ \circ \ \mathbf{v} \ + \ 
\Bigl(\mathbf{I} \otimes \mathbf{I}\Bigl) \ \circ \ \Bigl( {}^{(1, x)}\boldsymbol{\Gamma} \cdot \mathbf{v}\Bigl) \\ \\ 
\frac{\partial \mathbf{R}_\rho}{\partial \boldsymbol{v}} = \mathbf{K}_{\rho v} =& \ \Bigl({}^{(1, x)}\boldsymbol{\Gamma} \cdot \boldsymbol{\rho}\Bigl) \ \circ \ \Bigl(\mathbf{I} \otimes \mathbf{I}\Bigl) \ + \ \boldsymbol{\rho} \ \circ \ \Bigl({}^{(1, x)}\boldsymbol{\Gamma} \otimes \mathbf{I}\Bigl) \\ \\ 
\frac{\partial \mathbf{R}_\rho}{\partial \boldsymbol{e}} = \mathbf{K}_{\rho e} =& \ \mathbf{0} \\ \\
\frac{\partial \mathbf{R}_v}{\partial \boldsymbol{\rho}} = \mathbf{K}_{v\rho} =& \ \Bigl(\mathbf{I} \otimes \mathbf{I} \Bigl) \ \circ \ \Bigl[\Bigl(\mathbf{v}\cdot {}^{(1, t)}\boldsymbol{\Gamma}\Bigl) \ + \ \Bigl({}^{(1, x)}\boldsymbol{\Gamma} \cdot \mathbf{v}\Bigl) \ \circ \ \mathbf{v} \Bigl] \\ \\ 
\frac{\partial \mathbf{R}_v}{\partial \boldsymbol{v}} = \mathbf{K}_{vv} =& \ \boldsymbol{\rho} \ \circ \ \Bigl[ \mathbf{I} \otimes {}^{(1, t)}\boldsymbol{\Gamma} \ + \ \Bigl({}^{(1, x)}\boldsymbol{\Gamma}\otimes \mathbf{I}\Bigl) \ \circ \ \mathbf{v} \ + \ \Bigl({}^{(1, x)}\boldsymbol{\Gamma}\cdot \mathbf{v}\Bigl) \ \circ \ \Bigl(\mathbf{I}\otimes \mathbf{I}\Bigl) \Bigl] \ - \ \frac{4}{3}\mu \Bigl({}^{(2, x)}\boldsymbol{\Gamma}\otimes \mathbf{I}\Bigl) \\ \\ 
\frac{\partial \mathbf{R}_v}{\partial \boldsymbol{e}} = \mathbf{K}_{ve} =& \ \Bigl(\gamma-1\Bigl)\Bigl[\Bigl({}^{(1, x)}\boldsymbol{\Gamma} \cdot \boldsymbol{\rho} \Bigl) \ \circ \ \Bigl(\mathbf{I} \otimes \mathbf{I}\Bigl) \ + \ \boldsymbol{\rho} \ \circ \ \Bigl({}^{(1, x)}\boldsymbol{\Gamma} \otimes \mathbf{I}\Bigl)\Bigl] \\ \\ 
\frac{\partial \mathbf{R}_e}{\partial \boldsymbol{\rho}} = \mathbf{K}_{e \rho} =& \ \Bigl(\mathbf{I} \otimes \mathbf{I}\Bigl) \ \circ \ \Bigl[ \mathbf{e} \cdot {}^{(1, t)}\boldsymbol{\Gamma} \ + \ \mathbf{v} \ \circ \ \Bigl({}^{(1, x)}\boldsymbol{\Gamma} \cdot \mathbf{e}\Bigl)\Bigl] \ + \ \Bigl[\Bigl(\gamma-1\Bigl)\Bigl(\mathbf{I} \otimes \mathbf{I}\Bigl) \ \circ \ \mathbf{e}\Bigl] \ \circ \ \Bigl({}^{(1, x)}\boldsymbol{\Gamma} \cdot \mathbf{v}\Bigl) \\ \\ 
\frac{\partial \mathbf{R}_e}{\partial \mathbf{v}} = \mathbf{K}_{e v} =& \ \boldsymbol{\rho} \ \circ \ \Bigl(\mathbf{I} \otimes \mathbf{I}\Bigl) \ \circ \ \Bigl({}^{(1, x)}\boldsymbol{\Gamma}\cdot \mathbf{e}\Bigl) \ - \ \frac{4}{3}\mu \Bigl({}^{(1, x)}\boldsymbol{\Gamma}\otimes \mathbf{I}\Bigl) \ \circ \ \Bigl({}^{(1, x)}\boldsymbol{\Gamma}\cdot \mathbf{v}\Bigl) \ - \\ & \Bigl[ \frac{4}{3}\mu \Bigl({}^{(1, x)}\boldsymbol{\Gamma}\cdot \mathbf{v}\Bigl)\ - \ \Bigl(\gamma-1\Bigl)\boldsymbol{\rho} \ \circ \ \mathbf{e}\Bigl] \ \circ  \ \Bigl({}^{(1, x)}\boldsymbol{\Gamma} \otimes \mathbf{I}\Bigl) \\ \\ 
\frac{\partial \mathbf{R}_e}{\partial \mathbf{e}} = \mathbf{K}_{e e} =& \ \boldsymbol{\rho} \ \circ \ \Bigl[\Bigl(\mathbf{I}\otimes {}^{(1, t)}\boldsymbol{\Gamma}\Bigl) \ + \ \mathbf{v} \ \circ \ \Bigl({}^{(1, x)}\boldsymbol{\Gamma} \otimes \mathbf{I}\Bigl)\Bigl] \ + \ \Bigl[\Bigl(\gamma-1\Bigl) \boldsymbol{\rho} \ \circ \ \Bigl(\mathbf{I} \otimes \mathbf{I}\Bigl)\Bigl] \ \circ \ \Bigl({}^{(1, x)}\boldsymbol{\Gamma}\cdot \mathbf{v}\Bigl) \ - \ \frac{\kappa}{\mathrm{c}_v}\Bigl({}^{(2, x)}\boldsymbol{\Gamma} \otimes \mathbf{I}\Bigl)
\end{split}
\end{align*}
\end{appendices}

\end{document}